\documentclass[12pt, english]{article} 

\usepackage{color}
\usepackage[margin=1in]{geometry} 
\usepackage{graphicx} 
\usepackage{amsmath} 
\usepackage{amsfonts} 
\usepackage{amsthm} 
\usepackage{amssymb}
\usepackage{mathrsfs}
\usepackage{mathtools}
\usepackage{enumitem}
\usepackage{babel}
\usepackage{hyperref}
\newtheorem{theorem}{Theorem}[section]

\newtheorem{proposition}[theorem]{Proposition}

\newtheorem{definition}[theorem]{Definition}
\newtheorem{remark}[theorem]{Remark}

\DeclareMathAlphabet{\mathpzc}{OT1}{pzc}{m}{it}

\providecommand{\keywords}[1]{\textbf{\textit{Keywords: }}#1}
\title{Discretizations of Stochastic Evolution Equations in Variational Approach Driven by Jump-Diffusion %
} 
\author{Sima Mehri\footnotemark[1]\and 
Erfan Salavati\footnotemark[2] \and 
Bijan Z. Zangeneh\footnotemark[3] 
} 

\begin{document} 

\maketitle
\renewcommand{\thefootnote}{\fnsymbol{footnote}}
\footnotetext[1]{Department of Computer Science, University of Warwick, Coventry, England, United Kingdom}
\footnotetext[2]{Department of Mathematics and Computer Science, Amirkabir University of Technology (Tehran Polytechnic), Tehran, Iran}
\footnotetext[3]{Department of Mathematical Sciences, Sharif University of Technology, Tehran, Iran}
\renewcommand{\thefootnote}{\arabic{footnote}}
\maketitle 

\begin{abstract}
Stochastic evolution equations with compensated Poisson noise are considered in the variational approach with monotone and coercive coefficients. Here the Poisson noise is assumed to be time homogeneous with $\sigma$-finite intensity measure on a metric space. By using finite element methods and Galerkin approximations, some explicit and implicit discretizations for this equation are presented and their convergence is proved.  Polynomial growth condition and linear growth condition are assumed on the drift operator, respectively for the implicit and explicit schemes.

\end{abstract}
\section{Introduction}
\label{intro} 
There exist many results for numerical approximations of stochastic evolution equations (SEEs) driven by both continuous and c\`adl\`ag martingales. In particular, for stochastic partial differential equations (SPDEs) driven by Wiener noise, results exist concerning numerical approximations in variational setting (see \cite{breit2021space}, \cite{emmrich2017nonlinear}, \cite{gyongy2005solutions}, \cite{gyongy2005discretization}, \cite{gyongy2007rate}, \cite{gyongy2009rate}, \cite{gyongy2016convergence}, \cite{kruse2021bdf2}, \cite{sauer2015lattice}, \cite{sauer2016analysis}), approximations of linear stochastic partial differential equations (see e.g. \cite{gerencser2015finite}, \cite{gerencser2017localization}, \cite{hall2012accelerated} and references therein) and numerical approximation of SPDEs in the semigroup framework (see e.g. \cite{abdulle2021convergence}, \cite{adamu2011numerical}, \cite{brehier2020influence}, \cite{cox2013pathwise}, \cite{de2006weak}, \cite{debbi2011space}, \cite{debussche2011weak}, \cite{debussche2009weak}, \cite{gyongy1998lattice}, \cite{gyongy1999lattice}, \cite{jentzen2011higher}, \cite{jentzen2011efficient}, \cite{jentzen2009numerical}, \cite{jentzen2020strong}, \cite{kruse2013strong}, \cite{kurniawan2018weak}, \cite{lindgren2012weak}, \cite{lord2010stochastic}, \cite{milstein2009solving}, \cite{mishura2007approximation}, \cite{quer2006space}, \cite{salimova2019numerical}, \cite{tambue2016weak}, \cite{wagner2008optimal}, \cite{wrzosek2019newton}, \cite{yang2015full}). 

More recently results concerning explicit schemes for equations with drift operator not satisfying the linear growth condition have been obtained for finite dimensional stochastic differential equations (SDEs) driven by Wiener noise in \cite{beyn2016stochastic}, \cite{beyn2017stochastic}, \cite{halidias2014novel}, \cite{halidias2015construction}, \cite{halidias2016numerical}, \cite{hutzenthaler2015numerical}, \cite{hutzenthaler2020perturbation}, \cite{hutzenthaler2012strong}, \cite{hutzenthaler2018exponential}, \cite{kumar2019milstein}, \cite{liu2013strong}, \cite{mao2015truncated}, \cite{mao2016convergence}, \cite{ngo2017strong}, \cite{sabanis2013note}, \cite{sabanis2016euler}, \cite{song2018convergence}, \cite{szpruch2013v}, \cite{tretyakov2013fundamental}, \cite{wang2013tamed}, \cite{zhang2014new}, \cite{zong2014convergence} and for SDEs with jump-diffusive noise in \cite{dareiotis2016tamed}, \cite{kumar2014tamed}, \cite{tambue2015strong}. 

For semilinear SPDEs with super-linearly growing drift driven by Wiener or Gaussian noise, numerical schemes are considered in \cite{beccari2019strong}, \cite{bessaih2018numerical}, \cite{campbell2018adaptive}, \cite{feng2021strong}, \cite{hutzenthaler2020perturbation}, \cite{hutzenthaler2018strong}, \cite{jentzen2009pathwise}, \cite{jentzen2020strong}, \cite{jentzen2015milstein}, \cite{kamrani2016numerical}, \cite{kamrani2017pathwise}, \cite{kamrani2018implicit}, \cite{kloeden2011exponential}, \cite{liu2020strong}, \cite{printems2001discretization}, \cite{puvsnik2020strong}, \cite{yang2017convergence}, \cite{zhang2017numerical}.

Numerical schemes for SPDEs driven by Wiener/Gaussian noise, having super-linearly growing drift operator in variational setting are given in \cite{breit2021space}, \cite{emmrich2017nonlinear}, \cite{gyongy2005discretization}, \cite{gyongy2016convergence}, \cite{kruse2021bdf2}, \cite{sauer2015lattice}, \cite{sauer2016analysis}. 

For special types of SPDEs with super-linearly growing drift with Gaussian noise, numerical schemes have been obtained e.g. for Navier-Stokes equations in \cite{bessaih2018strong}, \cite{breckner2000galerkin}, \cite{brzezniak2013finite}, \cite{carelli2012numerical}, \cite{carelli2012rates}, \cite{dorsek2012semigroup}, \cite{duan2013finite}, \cite{li2018error}, \cite{mazzonetto2018strong}, 
for incompressible Euler equations in \cite{hong2021splitting}, 
for equations of geophysical fluid dynamics in \cite{glatt2017time}, for stochastic total variation flow in
\cite{rockner2021convergent}, for Ginzburg-Landau equations in \cite{becker2019strong}, for Allen-Cahn equations in
\cite{becker2017strong}, \cite{blomker2019numerically}, \cite{brehier2019strong}, \cite{brehier2019analysis}, \cite{kovacs2015backward}, \cite{kovacs2018discret}, \cite{ma2021convergence}, \cite{wang2020efficient}, 
for heat equation in \cite{becker2020lower}, 
for Cahn-Hilliard equations in \cite{furihata2018strong}, \cite{kossioris2012finite}, for Burgers equations in \cite{blomker2013galerkin}, \cite{blomker2013full}, \cite{ghayebi2017numerical}, \cite{hutzenthaler2019strong}, \cite{jentzen2019strong}, \cite{kamrani2012spectral}
, for Schr{\"o}dinger equations in \cite{cui2019strong}, and for Manakov equations in \cite{gazeau2013strong}, \cite{gazeau2014probability}.


For stochastic partial differential equations driven by c\`adl\`ag martingales in the semigroup framework one should mention e.g. \cite{barth2012milstein}, \cite{barth2019stochastic}, \cite{barth2018weak}, \cite{hausenblas2008finite}, \cite{kovacs2015weak}, \cite{lindner2013weak}, \cite{quer2005stochastic} (and references therein).
Numerical schemes for linear stochastic integro-differential equations of parabolic type arising in non-linear filtering of jump-diffusion processes have been obtained in \cite{dareiotis2016finite}. Numerical approximation for stochastic Schr{\"o}dinger type equations driven by martingale noise has been obtained in \cite{brzezniak2012splitting}.

To the best of our knowledge, numerical schemes for (general) stochastic partial differential equations with Poisson noise/L\'evy noise have not yet been considered in the variational setting.

In the article \cite{gyongy2005discretization}, I. Gy\"ongy and A. Millet studied discretizations of stochastic partial differential equations with Wiener noise in the variational setting. In this paper, we will generalize their approach by adding also compensated Poisson noise. More precisely, we consider the equation
\[ 
u_t=\zeta +\int_0^t{A_s\left(\bar{u}_s\right)\mathrm{d} s}+\int_0^t{B_s\left(\bar{u}_s\right)\mathrm{d} W_s}+\int_0^t\int_E{F_s\left(\bar{u}_{s},\xi\right)\tilde{N}\left(\mathrm{d} s,\mathrm{d}\xi\right)}, 
\] 
with respect to a Gelfand triple $V \hookrightarrow H \hookrightarrow V^*$. We assume $V$ to be a reflexive separable real Banach space which is continuously and densely embedded into the Hilbert space $H$. $V^*$ stands for the dual space of $V$ containing $H^*$, the dual of $H$, as a dense subset. Identifying $H$ with $H^*$, we obtain the dense and continuous embeddings $V \hookrightarrow H \hookrightarrow V^*$. A solution of the equation above is supposed to be a c\`adl\`ag $H$-valued stochastic process $u$ taking values in $V$ almost everywhere and therefore has a $V$-valued predictable modification $\bar{u}$. 

In the above equation $W$ is a cylindrical Wiener process in a Hilbert space $U$. We denote by ${N}$ a Poisson random measure, independent of $W$, which is considered to be time homogeneous with $\sigma$-finite intensity measure $\nu$ on a metric space $E$. $\tilde{N}$ is the compensated Poisson random measure corresponding to $N$. The coefficients $A$, $B$ and $F$ take values in $V^*$, $L_2(U,H)$ and $H$ respectively. Here, $L_2(U,H)$ is the space of Hilbert-Schmidt operators from $U$ to $H$. In this article, following \cite{gyongy2005discretization}, we assume hemicontinuity and growth condition on $A$ (cf. \ref{C4} and \ref{C3} below), and monotonicity and coercivity conditions on $(A,B,F)$ (cf. \ref{C1} and \ref{C2}). 

Following \cite{gyongy2005discretization}, we provide explicit and implicit numerical schemes for the above equation. Here, the condition for convergence of explicit scheme is weaker than the corresponding condition in \cite{gyongy2005discretization}. We will consider numerical schemes with respect to equipartition of the time interval $[0,T]$ into subintervals $[t_{i-1}, t_i]$ and the approximated value of $u$ at time $t_i$ will be calculated implicitly or explicitly. In the explicit scheme, the operators $A$, $B$ and $F$ at every time will be replaced with their integral means, taken on the previous time subinterval. Orthogonal projection to a finite dimensional subspace of $V$ with respect to the inner product of $H$ is also essential for this explicit scheme. Similarly, in the implicit scheme, the operators $B$ and $F$ will be replaced with their integral means taken on the preceding time subintervals, but $A$ will be replaced with its integral mean taken on the current time subinterval. The orthogonal projection to finite dimensional subspaces is optional in the implicit method and this will give us two types of implicit schemes.

For the mathematical background on the variational approach to stochastic evolution equations we refer to \cite{krylov1981stochastic}, \cite{liu2015stochastic}, \cite{pardoux1975equations}, and \cite{prevot2007concise}. For stochastic partial differential equations driven by L\'evy noise or more general Poisson random measures we refer to the books \cite{applebaum2009levy} and \cite{peszat2007stochastic}.

\section{About the Equation} 
\label{sect:basics} 
Let $V$ be a reflexive separable Banach space, embedded densely and continuously in a Hilbert space $H$. Its dual space $H^\ast$ is then densely and continuously embedded in $V^*$. Identifying $H$ with its dual space $H^\ast$ using Riesz' isometry, we obtain the Gelfand triple $V \hookrightarrow H \hookrightarrow V^*$. We denote by $\left\langle \, ,\right\rangle$ the duality between $V$ and $V^*$ and by $\left(\, ,\right)_H$ the inner product in $H$. 
We will be interested in equations of the following type 
\begin{equation} \label{equation1} 
	u_t=\zeta +\int_0^t{A_s\left(\bar{u}_s\right)\mathrm{d} s}+\int_0^t{B_s\left(\bar{u}_s\right)\mathrm{d} W_s}+\int_0^t\int_E{F_s\left(\bar{u}_s,\xi\right)\tilde{N}\left(\mathrm{d} s,\mathrm{d}\xi\right)},
\end{equation} 
where $u$ is a c\`adl\`ag $H$-valued process and $\bar{u}$ is its predictable $V$-valued modification. 

Let $\left( \Omega,\mathcal{F},\mathcal{F}_t, \mathbb{P}\right)$ be a complete probability space such that the filtration $\mathcal{F}_t$ satisfies the usual conditions, i.e. it is right continuous and $\mathcal{F}_0$ contains all $P$-null sets. 

Let $W$ be an adapted cylindrical Wiener process in a Hilbert space $U$ such that for $t>s$, $W_t-W_s$ is independent of $\mathcal{F}_s$. 
Let ${N}$ be an adapted time homogeneous Poisson random measure on 
\[
\left(\left[0,T\right] \times E, \mathcal{B}\left(\left[0,T\right]\right)\otimes \mathcal{E}\right),
\] 
independent of $W$. Here $E$ is a metric space and $\mathcal{E}$ is its Borel $\sigma$-field. 
We assume that the intensity measure $\nu$ of $N$ is $\sigma$-finite on the metric space $\left(E,\mathcal{E}\right)$ where $E$ is countable union of compact sets, and that ${N}\left(\left(s,t\right],.\right)$ is independent of $\mathcal{F}_s$, as for the Wiener process. Finally, let $\tilde{N} := 
N-\mathrm{d} t\otimes \nu$ be the compensated Poisson random measure associated with $N$. 

In the next step let us specify the measurability assumptions on the coefficients $A$, $B$, and $F$. Using the notation of \cite{liu2010spde}, let $\mathcal{BF}$ be the $\sigma$-field of progressively measurable sets on $\left[0,T\right]\times\Omega$, i.e. 
\[\mathcal{BF}=\left\lbrace A\subseteq\left[0,T\right]\times\Omega : \forall t\in\left[0,T\right],A\cap\left(\left[0,t\right]\times\Omega\right)\in\mathcal{B}\left(\left[0,t\right]\right)\otimes\mathcal{F}_t\right\rbrace.\] 
Let $\mathcal{P}$ denote the predictable $\sigma$-field, i.e. the $\sigma$-field generated by all left continuous and $\mathcal{F}_t$-adapted real-valued processes on $\left[0,T\right]\times\Omega$. We denote by $\left(L_{2}\left(U,H\right),\left\langle \, , \right\rangle_2,\left\lVert\, \right\rVert_2\right)$ the space of Hilbert-Schmidt operators from the Hilbert space $U$ to $H$. 
Then we assume that 
\[A:\left(\left[0,T\right]\times\Omega\times V,\mathcal{BF}\otimes\mathcal{B}\left(V\right)\right)\rightarrow\left(V^*,\mathcal{B}\left(V^*\right)\right),\] 
\[B:\left(\left[0,T\right]\times\Omega\times V,\mathcal{BF}\otimes\mathcal{B}\left(V\right)\right)\rightarrow\left(L_2(U,H),\mathcal{B}\left(L_2(U,H)\right)\right),\] 
\[F:\left(\left[0,T\right]\times\Omega\times V\times E,\mathcal{P}\otimes\mathcal{B}\left(V\right)\otimes\mathcal{E}\right)\rightarrow\left(H,\mathcal{B}\left(H\right)\right)\] 
are measurable. Note that predictability of $F$ is required, since it is integrated with respect to the Poisson random measure. 
Now we shall make some assumptions on the operators $A$, $B$ and $F$ and the initial condition $\zeta $. Let $p\in[2,\infty)$ and $q$ be its conjugate i.e. $1/p+1/q=1$. Let $K_1, K_2$ be non-negative integrable functions on $[0,T]\times \Omega$ and $\lambda$ and $\mu$ be respectively positive and non-negative deterministic integrable functions on $[0,T]$. The following conditions will be needed throughout the paper: 
\begin{enumerate}[label={(C\theenumi)}] 
	\item \label{C1} Monotonicity condition on $(A,B,F)$: almost surely for all $x,y\in V$ and all $ t\in[0,T]$, 
	\begin{align*} 
		\MoveEqLeft 2\left\langle A_t(x)-A_t(y),x-y\right\rangle+\| B_t(x)-B_t(y)\|^2_2 \\&+\int_E\|F_t(x,\xi)-F_t(y,\xi)\|_H^2\nu(\mathrm{d}\xi)\leq 0.
	\end{align*} 
	\item \label{C2} Coercivity condition on $(A,B,F)$: almost surely for all $x\in V$ and all $ t\in[0,T]$, 
	\begin{align*}\MoveEqLeft 2\left\langle A_t(x),x\right\rangle+\|B_t(x)\|^2_2+\int_E\|F_t(x,\xi)\|_H^2\nu(\mathrm{d}\xi)+\lambda(t)\|x\|_V^p \\ & \leq K_1(t)+\mu(t)\|x\|_H^2.
	\end{align*} 
	\item \label{C3} Growth condition on $A$: there exists $\alpha> 0$ such that almost surely for all $x\in V$ and all $t\in[0,T]$, 
	\[\|A_t(x)\|^q_{V^*}\leq\alpha\lambda^q(t)\|x\|_V^p+K_2(t)\lambda^{q-1}(t).\] 
	\item \label{C4} Hemicontinuity of $A$: almost surely for all $x,y,z\in V$ and all $t\in[0,T]$, 
	\[\lim_{\epsilon\rightarrow 0}\left\langle A_t(x+\epsilon y),z\right\rangle=\left\langle A_t(x),z\right\rangle.\] 
	\item\label{C6} $\zeta \in L^2(\Omega, \mathcal{F}_0;H)$. 
	\item\label{C7} There exists an increasing sequence 
	$E^1\subset E^2\subset E^3\subset \cdots$ of compact subsets of $E$, having finite $\nu$-measure, such that $\bigcup_{l=1}^\infty E^l=E$.
\end{enumerate} 

The monotonicity condition \ref{C1} can be weakened as follows:
\begin{enumerate}[label={({C}\theenumi ')}] 
	\item \label{Ctilde1} Almost surely for all $x,y\in V$ and all $ t\in[0,T]$, 
	\begin{align*} 
		2&\left\langle A_t(x)-A_t(y),x-y\right\rangle +\| B_t(x)-B_t(y)\|^2_2 
	\\ 
	&	+\int_E\|F_t(x,\xi)-F_t(y,\xi)\|_H^2\nu(\mathrm{d}\xi) 
		 \leq K(t)\left\lVert x-y\right\rVert^2_H,
	\end{align*}
	where $K$ is a non-negative integrable function on $[0,T]$. 
\end{enumerate}

Indeed, let $u_t$ be a solution to equation \eqref{equation1} and let $\gamma_t:=\exp{\left(-\frac{1}{2}\int_0^t K_s\mathrm{d} s\right)}$. Then $v_t=\gamma_t^{-1}u_t$ is a solution of 
\[
v_t=\zeta +\int_0^t{\bar{A}_s\left(\bar{v}_s\right)\mathrm{d} s}+\int_0^t{\bar{B}_s\left(\bar{v}_s\right)\mathrm{d} W_s}+\int_0^t\int_E{\bar{F}_s\left(\bar{v}_s,\xi\right)\tilde{N}\left(\mathrm{d} s,\mathrm{d}\xi\right)} 
\]
where 
\begin{align*} 
	\bar{A}_t(x):=\gamma_t^{-1}A_t(\gamma_t x)-\frac{1}{2}K_t x, \bar{B}_t(x):=\gamma_t^{-1}B_t(\gamma_t x), \bar{F}_t(x,\xi):=\gamma_t^{-1}F_t(\gamma_t x,\xi). 
\end{align*}
If $A$, $B$, and $F$ satisfy \ref{Ctilde1}, \ref{C2}, \ref{C4} it follows that $\bar{A},\bar{B}$ and $\bar{F}$ satisfy \ref{C1}, \ref{C2}, \ref{C4}. If $A$ satisfies \ref{C3} and $K(t)\leq C\lambda(t), t\in [0,T]$ for some constant $C$, then $\bar{A}$ satisfies \ref{C3}.

\begin{proposition}\label{proposition1} 
	Condition \ref{C3} together with each of \ref{C1} and \ref{C2} gives respectively the following conditions on $(B,F)$: 
	\begin{equation}\label{propBF2}
		\begin{split}
			\MoveEqLeft\| B_t(x)-B_t(y)\|^2_2 +\int_E\|F_t(x,\xi)-F_t(y,\xi)\|_H^2\nu(\mathrm{d}\xi)\\&\leq 3\left(\alpha^{1/q} +\frac{1}{p}\right)\lambda(t)\left(\left\lVert x\right\rVert_V^p+\left\lVert y\right\rVert_V^p\right)+\frac{8}{q} K_2(t),
		\end{split}
	\end{equation}
	
	\begin{equation}\label{propBF1}
		\begin{split}\MoveEqLeft\left\lVert B_t(x)\right\rVert^2_2+\int_E\left\lVert F_t(x,\xi)\right\rVert_H^2\nu(\mathrm{d}\xi)\\&\leq 6\left(\alpha^{1/q} +\frac{1}{p}\right)\lambda(t)\left\lVert x\right\rVert_V^p+\frac{16}{q} K_2(t) +2K_1(t).
		\end{split}
	\end{equation}
\end{proposition} 
\begin{proof} 
	We get by \ref{C3} that 
	\[\left\lVert A_t(x)\right\rVert_{V^*}\leq\alpha^{1/q}\lambda(t)\left\lVert x\right\rVert_V^{p-1}+K_2^{1/q}(t)\lambda^{1/p}(t).\] 
	So using Young's inequality we have 
	\begin{align*} 
		\left\lvert\left\langle x,A_t(x) \right\rangle\right\rvert &\leq\alpha^{1/q}\lambda(t)\left\lVert x\right\rVert_V^p+K_2^{1/q}(t)\lambda^{1/p}(t)\left\lVert x\right\rVert_V\\&\leq\left(\alpha^{1/q}+\frac{1}{p}\right)\lambda(t)\left\lVert x\right\rVert_V^p+\frac{1}{q}K_2(t)
	\end{align*} 
	and 
	\begin{align*} 
		\left\lvert\left\langle y,A_t(x) \right\rangle\right\rvert &\leq\alpha^{1/q}\lambda(t)\left\lVert x\right\rVert_V^{p-1}\left\lVert y\right\rVert_V+K_2^{1/q}(t)\lambda^{1/p}\left\lVert y\right\rVert_V\\&\leq \alpha^{1/q}\lambda(t)\left( \frac{p-1}{p}\left\lVert x\right\rVert_V^p+\frac{1}{p}\left\lVert y\right\rVert_V^p\right)+
		\frac{1}{p}\lambda(t)\left\lVert y\right\rVert_V^p+\frac{1}{q}K_2(t).
	\end{align*}
	Combining these inequalities with \ref{C1} and \ref{C2} yields
	\begin{equation*}
		\begin{split}
			\MoveEqLeft\| B_t(x)-B_t(y)\|^2_2 +\int_E\|F_t(x,\xi)-F_t(y,\xi)\|_H^2\nu(\mathrm{d}\xi)\\&\leq 2\left\lvert\left\langle x, A_t(y)\right\rangle\right\rvert+2\left\lvert\left\langle y, A_t(x)\right\rangle\right\rvert+ 2\left\lvert\left\langle x, A_t(x)\right\rangle\right\rvert+2\left\lvert\left\langle y, A_t(y)\right\rangle \right\rvert\\&\leq 
				3\left(\alpha^{1/q}+\frac{1}{p}\right)\lambda(t)\left( \left\lVert x\right\rVert_V^p+\left\lVert y\right\rVert_V^p\right)+\frac{8}{q}K_2(t)
		\end{split}
	\end{equation*}
	and
		\begin{align*} 
			\MoveEqLeft\left\lVert B_t(x)\right\rVert^2_2+\int_E\left\lVert F_t(x,\xi)\right\rVert_H^2\nu(\mathrm{d}\xi)\\&\leq 2\| B_t(x)-B_t(0)\|^2_2 +2\int_E\|F_t(x,\xi)-F_t(0,\xi)\|_H^2\nu(\mathrm{d}\xi)\\&\quad + 2\left\lVert B_t(0)\right\rVert^2_2+2\int_E\left\lVert F_t(0,\xi)\right\rVert_H^2\nu(\mathrm{d}\xi)\\&\leq 6\left(\alpha^{1/q} +\frac{1}{p}\right)\lambda(t)\left\lVert x\right\rVert_V^p+\frac{16}{q} K_2(t) +2K_1(t).
	\end{align*}
\end{proof} 
Now we are going to define the solution of equation \eqref{equation1}. First we remind the notion of modification of a stochastic process. 
\begin{definition} 
	Let $z$ be a stochastic process. $\bar{z}$ is called a modification of $z$ if for $\mathrm{d} t\otimes\mathbb{P}$-almost all $(t,\omega)\in [0,T]\times\Omega$, $\bar{z}(t,\omega)=z(t,\omega)$. 
\end{definition}

\begin{definition}\label{solutiondef} 
	A c\`adl\`ag $H$-valued $(\mathcal{F}_t)$-adapted stochastic process $u$ is a (strong) solution to the equation \eqref{equation1}, if it has a predictable $V$-valued modification $\bar{u}$ in 
	\[L^p([0,T]\times\Omega, \lambda\mathrm{d} t \otimes \mathbb{P};V)\cap L^2([0,T]\times\Omega, \mathrm{d} t \otimes \mathbb{P};H)\] 
	such that the equation 
	\begin{align*}\left( u_t,v\right)_H&=\left( \zeta ,v\right)_H+\int_0^t\left\langle A_s(\bar{u}_s),v\right\rangle \mathrm{d} s 
		+\int_0^t\left( v, B_s(\bar{u}_s) \mathrm{d} W_s\right)_H\\&\quad+\int_0^t\int_E\left( F_s(\bar{u}_s,\xi),v\right)_H\tilde{N}(\mathrm{d} s,\mathrm{d}\xi) 
	\end{align*} 
	holds for all $v \in V$ and $\mathrm{d} t\otimes \mathbb{P}$- almost every $(t,\omega) \in [0,T]\times\Omega$. 
\end{definition} 

\begin{remark} 
	Suppose that $z$ is an adapted c\`adl\`ag stochastic process in H that $\mathrm{d} t\otimes\mathbb{P}$-almost everywhere belongs to $V$. Since $V$ is a Borel subset of $H$ and $z(t-), t\in [0,T]$ is a predictable modification of $z$ in $H$, we have $z(t-)\mathbf{1}_{\left\lbrace z(t-)\in V\right\rbrace}$ as a $V$-valued predictable modification of $z$.
\end{remark} 
The following existence and uniqueness theorems hold (see \cite{brzezniak2014strong}, \cite{gyongy1982stochastic}).
\begin{theorem}\label{exist} 
	Let conditions \ref{C1}-\ref{C7} hold. Then equation \eqref{equation1} has a solution $u$ that satisfies
	\begin{equation}\label{Boundedness of Solution} 
		\sup_{t\in[0,T]}\mathbb{E}\left\lVert u_t\right\rVert_H^2<\infty. 
	\end{equation}
\end{theorem}

\begin{theorem}\label{unique}
	Assuming \ref{C1}-\ref{C7}, the solution of \eqref{equation1} in the sense of Definition \ref{solutiondef} is unique and satisfies  almost surely, for all $t\in[0,T]$,
	\begin{equation}\label{strong} 
		u_t=\zeta +\int_0^t A_s(\bar{u}_s)\mathrm{d} s+\int_0^t B_s(\bar{u}_s)\mathrm{d} W_s+\int_0^t\int_E F_s(\bar{u}_s,\xi)\tilde{N}(\mathrm{d} s,\mathrm{d} \xi) .
	\end{equation} 
\end{theorem} 

We will prove Theorem \ref{exist} by numerical approximations. To prove Theorem \ref{unique}, we use the following theorem of \cite{gyongy1982stochastics} with the constant stopping time $\tau\equiv T$. 

\begin{theorem}\label{theorem1of1} 
	\cite[Theorem 1]{gyongy1982stochastics} Let $\Lambda$ be an increasing adapted real valued stochastic process with c\`adl\`ag trajectories. Assume $z$ and $y$ are respectively $V$ and $V^*$-valued progressively measurable stochastic processes. Suppose that $\left\lVert z(t)\right\rVert_V$,$\left\lVert y(t)\right\rVert_{V^*}$ and $\left\lVert z(t)\right\rVert_V\times\left\lVert y(t)\right\rVert_{V^*}$ are locally integrable with respect to $\mathrm{d}\Lambda_t$, i.e. their trajectories are almost surely integrable with respect to $\mathrm{d}\Lambda_t$. Let $h(t)$ be an $H$-valued locally square integrable c\`adl\`ag martingale and $\tau$ denote a stopping time. Set $D=\left\lbrace(t,\omega)\in [0,T]\times\Omega:t\leq\tau(\omega)\right\rbrace$ and suppose that for all $v\in V$, and $\mathrm{d}\Lambda_t\otimes \mathbb{P} $- almost every $(t,\omega) \in D$ we have 
	\begin{equation}\label{equ2} 
		\left(v,z(t)\right)_H=\int_0^t\left\langle v,y(s)\right\rangle \mathrm{d}\Lambda_s+\left(v,h(t)\right)_H. 
	\end{equation} 
	Then there exists a subset $\Omega'\subseteq\Omega$ with $\mathbb{P}(\Omega')=1$ and an adapted c\`adl\`ag $H$-valued process like $\tilde{z}$ that is equal to z for $\mathrm{d} \Lambda \otimes \mathbb{P}$ almost all $(t,\omega)$ and has the following property instead of \eqref{equ2}: 
	\begin{equation} 
		\forall\omega\in\Omega', \forall t \leq \tau(\omega), \forall v\in V: \\ 
		\left( v, \tilde{z}(t)\right)_H = \int_0^t\left\langle v, y(s)\right\rangle \mathrm{d}\Lambda_s+\left(v,h(t)\right)_H,
	\end{equation} 
	and the following It\^o formula for $\left\lVert\cdot\right\rVert_H^2$ holds: 
	\begin{align*}\left\lVert\tilde{z}(t)\right\rVert_H^2 &= \left\lVert h(0)\right\rVert_H^2+2\int_0^t\left\langle z(s), y(s)\right\rangle \mathrm{d}\Lambda_s+2\int_0^t\left(\tilde{z}(s-),\mathrm{d} h(s)\right)_H\\&\quad-\int_0^t\left\lVert y(s)\right\rVert_H^2\Delta \Lambda(s) \mathrm{d}\Lambda(s)+\left[h\right]_t ,
	\end{align*} 
	where $\left[h\right]_t$ is the quadratic variation of $h$ and $\Delta \Lambda(s)$ is the value $\Lambda(s)-\Lambda(s-)$. If $y(s)\notin H$, we set $\left\lVert y(s)\right\rVert_H^2:=\infty$. 
\end{theorem}

\begin{proof}[Proof of Theorem \ref{unique}]
	We first prove the uniqueness of the solution. Let $u^{(1)}$ and $u^{(2)}$ be two solutions for equation \eqref{equation1} with predictable $V$-valued modifications $\bar{u}^{(1)}$ and $\bar{u}^{(2)}$ respectively. Let us calculate $\left\lVert{u}_t^{(1)}-{u}_t^{(2)}\right\rVert^2_H$ by means of the It\^o formula from the previous theorem. For this purpose, we set $z(t) = \bar{u}^{(1)}_t-\bar{u}^{(2)}_t$ and define $\mathrm{d}\Lambda_t:=\lambda(t)\mathrm{d} t$. We have
	\[\left(v,z(t)\right)_H=\int_0^t\left\langle v,y(s)\right\rangle \mathrm{d}\Lambda_s+\left(v,h(t)\right)_H\quad \mathrm{d} t\otimes \mathbb{P}\text{-a.e.}\] 
	where 
	\begin{align*}
		y(t)&=\left(A_t\left(\bar{u}_t^{(1)}\right)-A_t\left(\bar{u}_t^{(2)}\right)\right)\lambda(t)^{-1},\\h(t)&=\int_0^t \left(B_s\left(\bar{u}_s^{(1)}\right)-B_s\left(\bar{u}_s^{(2)}\right)\right) \mathrm{d} W_s\\&\quad+\int_0^t\int_E\left({ F_s\left(\bar{u}_s^{(1)},\xi\right)}-{F_s\left(\bar{u}_s^{(2)},\xi\right)}\right)\tilde{N}(\mathrm{d} s,\mathrm{d}\xi). 
	\end{align*} 
	Note that $\mathrm{d} \Lambda_t\otimes \mathbb{P}$ and $\mathrm{d} t\otimes \mathbb{P}$ are absolutely continuous with respect to each other and therefore $\mathrm{d} \Lambda_t\otimes \mathbb{P}$-a.e. is equivalent to $\mathrm{d} t\otimes \mathbb{P}$-a.e. By \eqref{propBF2} we obtain that
	$h$ is a c\`adl\`ag square integrable $H$-valued martingale. It is obvious that $y$ and $z$ are progressively measurable. From the assumptions of Theorem \ref{theorem1of1}, it remains to check the local integrability of $\left\lVert z(t)\right\rVert_V$,$\left\lVert y(t)\right\rVert_{V^*}$ and $\left\lVert z(t)\right\rVert_V\left\lVert y(t)\right\rVert_{V^*}$ with respect to $\mathrm{d}\Lambda_t$. We have 
	\begin{align*} 
		\int_0^T\left\lVert y(t)\right\rVert_{V^*}^q \mathrm{d}\Lambda_t&=\int_0^T\left\lVert A_t\left(u_t^{(1)}\right)-A_t\left(u_t^{(2)}\right)\right\rVert_{V^*}^q\lambda^{1-q}(t)\mathrm{d} t\\&\leq 2^{q-1}\int_0^T\sum_{i=1}^2\alpha\lambda(t)\left\lVert u^{(i)}_t\right\rVert_V^p \mathrm{d} t+2^q \int_0^T K_2(t) \mathrm{d} t<\infty\quad \text{a.s.},\\\int_0^T\left\lVert z(t)\right\rVert_V^p \mathrm{d}\Lambda_t&=\int_0^T\left\lVert u_t^{(1)}-u_t^{(2)}\right\rVert_V^p\lambda(t)\mathrm{d} t\\&\leq 2^{p-1}\int_0^T\sum_{i-1}^2\left\lVert u^{(i)}_t\right\rVert_V^p\lambda(t)\mathrm{d} t<\infty\quad \text{a.s.}
	\end{align*} 
	Applying H\"older's inequality, we get 
	\[\begin{split}\MoveEqLeft[2]\int_0^T\left\lVert y(t)\right\rVert_{V^*}\left\lVert z(t)\right\rVert_V\mathrm{d}\Lambda_t\\&\leq\left(\int_0^T\left\lVert y(t)\right\rVert^q\mathrm{d}\Lambda_t\right)^{1/q}\left(\int_0^T\left\lVert z(t)\right\rVert^p\mathrm{d}\Lambda_t\right)^{1/p}<\infty\quad \text{a.s.}\end{split}\] 
	So the conditions of Theorem \ref{theorem1of1} with $\tau =T$ hold and ${u}^{(1)}-{u}^{(2)}$ which is the c\`adl\`ag $H$-valued modification of $z$ satisfies 
	\begin{equation} \label{mt}
		\begin{split}
			\left\lVert{u}_t^{(1)}-{u}_t^{(2)}\right\rVert^2_H &= 
			\int_0^t2\left\langle \bar{u}_s^{(1)}-\bar{u}_s^{(2)},A_s\left(\bar{u}_s^{(1)}\right)-A_s\left(\bar{u}_s^{(2)}\right)\right\rangle \mathrm{d} s\\&\quad +\int_0^t \left\lVert B_s\left(\bar{u}_s^{(1)}\right)-B_s\left(\bar{u}_s^{(2)}\right)\right\rVert_2^2\mathrm{d} s \\ 
			& \quad +\int_0^t\int_E\left\lVert{F_s\left(\bar{u}_s^{(1)},\xi\right)}-{F_s\left(\bar{u}_s^{(2)},\xi\right)}\right\rVert_H^2 \nu(\mathrm{d} \xi)\mathrm{d} s+m_t.
		\end{split} 
	\end{equation}
	where 
	\begin{align*} 
		m_t&=2\int_0^t\left({u}_{s-}^{(1)}-{u}_{s-}^{(2)},B_s\left(\bar{u}_s^{(1)}\right)-B_s\left(\bar{u}_s^{(2)}\right)\mathrm{d} W_s\right)_H\\ 
		&\quad+ 2\int_0^t\int_E \left({u}_{s-}^{(1)}-{u}_{s-}^{(2)},{F_s\left(\bar{u}_s^{(1)},\xi\right)}-{F_s\left(\bar{u}_s^{(2)},\xi\right)}\right)_H\tilde{N}(\mathrm{d} s,\mathrm{d}\xi)\\&\quad 
		+\int_0^t\int_E\left\lVert{F_s\left(\bar{u}_s^{(1)},\xi\right)}-{F_s\left(\bar{u}_s^{(2)},\xi\right)}\right\rVert_H^2\tilde{N}(\mathrm{d} s,\mathrm{d}\xi)
	\end{align*} 
	is a local martingale. The monotonicity condition \ref{C1} gives 
	\[\left\lVert{u}_t^{(1)}-{u}_t^{(2)}\right\rVert^2_H\leq m_t.\] 
	Let $\sigma_n\uparrow \infty$ be stopping times such that $m_{t\wedge \sigma_n},t\geq 0$ are martingales. Using Fatou's lemma, we have
	\[\mathbb{E} \left\lVert{u}_t^{(1)}-{u}_t^{(2)}\right\rVert^2_H\leq \liminf_{n\to\infty}\mathbb{E} \left\lVert{u}_{t\wedge \sigma_n}^{(1)}-{u}_{t\wedge \sigma_n}^{(2)}\right\rVert^2\leq \liminf_{n\to\infty} \mathbb{E}(m_{t\wedge \sigma_n})=0\]
	and therefore the uniqueness is proved. 
	
	 Let $u$ be a solution to equation \eqref{equation1}. In order to show equation \eqref{strong} for solution $u$, we apply Theorem \ref{theorem1of1} for $z=\bar{u}$, $\mathrm{d} \Lambda_t=\lambda(t)\mathrm{d} t$, $y=A_t(\bar{u}_t)\lambda(t)^{-1}$, $\tau=T$ and
		\[h(t)=\zeta +\int_0^tB_s\left(\bar{u}_s\right)\mathrm{d} W_s+\int_0^t\int_E F_s\left(\bar{u}_s,\xi\right)\tilde{N}(\mathrm{d} s,\mathrm{d} \xi).\] 
		By definition \ref{solutiondef}, equation \eqref{equ2} holds. By \eqref{propBF1}, we obtain that $h$ is a c\`adl\`ag square integrable $H$-valued martingale. It can be shown that $y$ and $z$ are progressively measurable. From the assumptions of Theorem \ref{theorem1of1}, it remains to check the local integrability of $\left\lVert z(t)\right\rVert_V$,$\left\lVert y(t)\right\rVert_{V^*}$ and $\left\lVert z(t)\right\rVert_V\left\lVert y(t)\right\rVert_{V^*}$ with respect to $\mathrm{d}\Lambda_t$. We have 
		\begin{align*} 
			\int_0^T\left\lVert y(t)\right\rVert_{V^*}^q \mathrm{d}\Lambda_t&=\int_0^T\left\lVert A_t(\bar{u}_t)\right\rVert_{V^*}^q\lambda^{1-q}(t)\mathrm{d} t\\&\leq \int_0^T\alpha\lambda(t)\left\lVert \bar{u}_t\right\rVert_V^p \mathrm{d} t+ \int_0^T K_2(t) \mathrm{d} t<\infty\quad \text{a.s.},\\\int_0^T\left\lVert z(t)\right\rVert_V^p \mathrm{d}\Lambda_t&=\int_0^T\left\lVert \bar{u}_t\right\rVert_V^p\lambda(t)\mathrm{d} t<\infty\quad \text{a.s.}
		\end{align*} 
		Applying H\"older's inequality, we get 
		\[\begin{split}\MoveEqLeft[2]\int_0^T\left\lVert y(t)\right\rVert_{V^*}\left\lVert z(t)\right\rVert_V\mathrm{d}\Lambda_t\\&\leq\left(\int_0^T\left\lVert y(t)\right\rVert^q\mathrm{d}\Lambda_t\right)^{1/q}\left(\int_0^T\left\lVert z(t)\right\rVert^p\mathrm{d}\Lambda_t\right)^{1/p}<\infty\quad \text{a.s.}\end{split}\] 
		So the conditions of Theorem \ref{theorem1of1} with $\tau =T$ hold. Therefore, we get almost surely, for all $t\in [0,T]$ and for all $v\in V$,
		\[\begin{split}(v,&u_t)_H=\int_0^t \left\langle v, A_s(\bar{u}_s)\right\rangle\mathrm{d} s+(v,\zeta)_H+\left(v,\int_0^t B_s(\bar{u}_s)\mathrm{d} W_s\right)_H\\&\quad\qquad+\left(v,\int_0^t\int_E F_s(\bar{u}_s,\xi)\tilde{N}(\mathrm{d} s,\mathrm{d} \xi)\right)_H\\&= \left(v,\zeta+\int_0^t A_s(\bar{u}_s)\mathrm{d} s+\int_0^t B_s(\bar{u}_s)\mathrm{d} W_s+\int_0^t\int_E F_s(\bar{u}_s,\xi)\tilde{N}(\mathrm{d} s,\mathrm{d} \xi)\right)_H\end{split}\]
		By using the fact that the space $V$ is a dense subset of $H$, this implies equation 
		\eqref{strong}.
\end{proof}

Suppose that $X$ is a separable Banach space, $\varphi$ is a positive function on $[0,T]$, and $p\in[1,\infty)$. For simplicity let us denote with $\mathscr{L}_X^p(\varphi)$ the space
\[L^p\left( [0,T]\times \Omega,\mathcal{B}([0,T])\otimes\mathcal{F}, \varphi\mathrm{d} t\otimes \mathbb{P};X\right).\]
Denote with $\mathscr{L}_X^p(\varphi,\mathcal{BF})$ and $\mathscr{L}_X^p(\varphi,\mathcal{P})$, the subspaces of $\mathscr{L}_X^p(\varphi)$ consisting of respectively progressively measurable and predictable processes. 
When $\varphi\equiv 1$, we use the notations $\mathscr{L}^p_X$, $\mathscr{L}^p_X(\mathcal{BF})$ and $\mathscr{L}^p_X(\mathcal{P})$. 
Denote by $\mathscr{G}$ the following Banach space 
\[\mathscr{G}:=\left\lbrace y\in \mathscr{L}_V^p(\lambda):\mathrm{ess}\sup_{t\in[0,T]}\mathbb{E}\left\lVert y_t\right\rVert^2_H<\infty\right\rbrace\] 
with the norm 
\[|y|_{\mathscr{G}}=\left(\mathbb{E}\int_0^T\left\lVert y_t\right\rVert^p_V\lambda(t)\mathrm{d} t\right)^{1/p}+\left(\mathrm{ess}\sup_{t\in[0,T]}\mathbb{E}\left\lVert y_t\right\rVert^2_H\right)^{1/2}.\] 
Let $\mathscr{G}_{\mathcal{BF}}$ be subspace of $\mathscr{G}$, consisting of progressively measurable processes. $\mathscr{G}_{\mathcal{BF}}$ is a Banach space too. Note that since for every $X$-valued adapted stochastic process like $z$, there exists a sequence of bounded continuous stochastic processes that converges to $z$ in $\mathscr{L}_X^p(\varphi)$, so $\mathscr{G}_{\mathcal{BF}}$ is dense in $\mathscr{L}_V^p(\lambda,\mathcal{BF})$.

Following \cite{gyongy2005discretization}, we characterize the solution of equation \eqref{equation1}. 

\begin{definition}\label{setA} 
	Denote by $\mathcal{A}$, the set consisting of quadruples $(\zeta,a,b,f)$ with the following conditions 
	\begin{enumerate}[label=(\roman{*})] 
		\item \label{setA:i}$\zeta\in L^2\left(\Omega,\mathcal{F}_0,\mathbb{P};H\right),$
		\item \label{setA:ii} $a \in \mathscr{L}^q_{V^*}(\lambda^{1-q},\mathcal{BF}),$ 
		\item \label{setA:iii}$b \in \mathscr{L}^2_{L_2(U,H)}(\mathcal{BF}),$ 
		\item \label{setA:iv}$f:\left([0,T]\times\Omega\times E,\mathcal{P}\otimes\mathcal{E}\right)\rightarrow\left(H,\mathcal{B}\left(H\right)\right)$ 
		satisfies 
		\[\mathbb{E}\int_0^T\int_E\left\lVert f(s,\xi)\right\rVert^2_H\nu(\mathrm{d}\xi)\mathrm{d} s<\infty,\] 
		\item \label{setA:v}there exists $x\in \mathscr{G}_{\mathcal{BF}}$ such that for all $v\in V$ and almost all $(t,\omega)\in [0,T]\times \Omega$, 
		\begin{align*} 
			(x_t,v)_H&=(\zeta,v)_H+\int_0^t \langle a_s,v\rangle \mathrm{d} s+\int_0^t \left(b_s\mathrm{d} W_s,v\right)_H\\&\quad+\int_0^t\int_E(f(s,\xi),v)_H\tilde{N}(\mathrm{d} s,\mathrm{d}\xi).
		\end{align*} 
	\end{enumerate} 
\end{definition} 
Let $(\zeta,a,b,f)$ belong to $\mathcal{A}$ and $x$ be the stochastic process as in part \ref{setA:v} of the above definition corresponding to the quadruple $(\zeta,a,b,f)$. For $y\in \mathscr{G}$, set 
\begin{align*} 
	I_y(\zeta,a,b,f)& := \mathbb{E}\int_0^T\Big[2\left\langle a_s-A_s(y_s),x_s-y_s\right\rangle+\left\lVert b_s-B_s(y_s)\right\rVert^2_2\\ 
	&\quad+\int_E\left\lVert f(s,\xi)-F_s(y_s,\xi) \right\rVert^2_H\nu(\mathrm{d}\xi)\Big] \mathrm{d} s\nonumber 
\end{align*} 
The next theorem which is an analogue of \cite[Theorem 2.7]{gyongy2005discretization} characterizes the solution of equation \eqref{equation1} and will be used for the proofs of the approximation theorems. 
\begin{theorem}\label{characterize} 
	Assume conditions \ref{C1}-\ref{C6}. If for some $(\zeta,a,b,f)\in\mathcal{A}$, and every $y\in \mathscr{G}_\mathcal{BF}$ \[I_y(\zeta,a,b,f)\leq 0,\] 
	then the stochastic process $x$ corresponding to $(\zeta,a,b,f)$ as in part \ref{setA:v} of Definition \ref{setA}, has an $H$-valued c\`adl\`ag modification which is a solution of equation \eqref{equation1} with initial condition $\zeta $ in the sense of Definition \ref{solutiondef}. 
\end{theorem} 
\begin{proof} 
	By Theorem \ref{theorem1of1}, $x$ has an adapted c\`adl\`ag $H$-valued modification, so it is sufficient to prove $a_s=A_s(x_s)$, $b_s=B_s(x_s)$ and $f(s,\xi)=F_s(x_s,\xi)$ almost everywhere. Set $y=x+\epsilon z$ where $z\in \mathscr{G}_\mathcal{BF}$, so $y\in\mathscr{G}_\mathcal{BF}$. Then
	\begin{equation}\label{ineq}
		\begin{split}
			0\geq I_y(\zeta,a,b,f)&  =  \mathbb{E}\int_0^T\Bigl[ 2\left\langle a_s-A_s(x_s+\epsilon z_s),-\epsilon z_s\right\rangle+\left\lVert b_s-B_s(x_s+\epsilon z_s)\right\rVert^2_2
			\\&\qquad+\int_E\left\lVert F_s(x_s+\epsilon z_s,\xi)-f(s,\xi)\right\rVert^2_H\nu(\mathrm{d}\xi)\Bigr] \mathrm{d} s
		\end{split} 
	\end{equation} 
	By choosing $\epsilon=0$, it follows that $b_s=B_s(x_s)$ and $f(s,\xi)=F_s(x_s,\xi)$ almost everywhere. Then \eqref{ineq} yields
	\[\mathbb{E}\int_0^T\left\langle a_s-A_s(x_s+\epsilon z_s),-z_s\right\rangle \mathrm{d} s\leq 0\]
	for all $z\in\mathscr{G}_\mathcal{BF}$. Hemicontinuity of $A$ implies 
	\[\lim_{\epsilon\rightarrow 0}\left\langle A_s(x_s+\epsilon z_s),z_s\right\rangle=\left\langle A_s(x_s),z_s\right\rangle,\] 
	and growth condition \ref{C3} implies that the functions $\left\langle A_s(x_s+\epsilon z_s),z_s\right\rangle$ are dominated by an integrable function on $[0,T]\times\Omega$. Hence, by dominated convergence,
	\[
	\mathbb{E}\int_0^T \left\langle a_s-A_s(x_s),-z_s\right\rangle \mathrm{d} s\leq 0
	\] 
	for every $z\in \mathscr{G}_\mathcal{BF}$. Substitute $-z$ instead of $z$, it is obvious that 
	\[\mathbb{E}\int_0^T \left\langle a_s-A_s(x_s),z_s\right\rangle \mathrm{d} s=0.\] 
	Since $a_s-A_s(x_s)$ belongs to $\mathscr{L}_{V^*}^q(\lambda^{1-q},\mathcal{BF})$ 
	and $\mathscr{G}_\mathcal{BF}$ is dense in $\mathscr{L}_V^p(\lambda,\mathcal{BF})$, the dual space of $\mathscr{L}_{V^*}^q(\lambda^{1-q},\mathcal{BF})$, we get $a_t=A_t(x_t)$ for almost all $(t,\omega)\in[0,T]\times\Omega$. 
\end{proof}

Now we are going to discretize space and time and the $\sigma$-finite measure $\nu$. Then we will apply these discretizations to the equation \eqref{equation1} and formulate explicit and implicit numerical schemes in the next section.

\section{Discretizations} 
Let us first introduce our space discretization. Let $V_1\subseteq V_2\subseteq\cdots$ be an increasing sequence of finite dimensional subspaces of $V$ such that $\bigcup_{n=1}^{\infty}V_n$ is dense in $V$. Consider the orthogonal projection operator $\Pi_n$ from $H$ onto $V_n$. Extend its domain to the space $V^*$ such that the operator remains continuous and linear and denote the obtained operator again by $\Pi_n$. Let \[\mathpzc{B}_n=\left\lbrace e_1,e_2,\ldots,e_{l_n}\right\rbrace\]
be a basis of $V_n$, orthonormal in $H$. 
 Then $\Pi_n$ has the following form: 
\[\forall x\in V^*\quad \Pi_nx=\left\langle e_1,x\right\rangle e_1+\left\langle e_2,x\right\rangle e_2+\cdots+\left\langle e_{l_n},x\right\rangle e_{l_n}\] 
\begin{proposition} 
	The following properties hold for $\Pi_n$: 
	\begin{enumerate}[label=(\roman{*})] 
		\item \label{Pi_n(i)}$\forall x \in V_n, \quad \Pi_n x=x$,
		\item \label{Pi_n(ii)}$\forall x \in V, y\in V^*,\quad \left\langle\Pi_n x,y\right\rangle=\left\langle x, \Pi_n y\right\rangle$,
		\item \label{Pi_n(iii)}$\forall h,k \in H,\quad \left(\Pi_n h,k\right)_H=\left( h,\Pi_nk\right)_H$, 
		\item \label{Pi_n(iv)}$\forall h\in H, \quad \left\lVert\Pi_n h\right\rVert_H\leq\left\lVert h\right\rVert_H, \lim_{n\rightarrow \infty}\left\lVert h-\Pi_n h\right\rVert_H=0$.
	\end{enumerate} 
\end{proposition} 

Let $\left\lbrace g_1,g_2,\ldots\right\rbrace$ be an orthonormal basis of $U$ and $\tilde{\Pi}_l$ be the orthogonal projection from $U$ to $\mathrm{span}\left\lbrace g_1,g_2,\ldots,g_l\right\rbrace$. Set 
\[W^l_t:=\tilde{\Pi}_lW_t=\sum_{k=1}^l(W_t, g_k)_U g_k.\] 
To approximate the compensated Poisson random measure, we use assumption \ref{C7}, which says there exists an increasing sequence 
$E^1\subset E^2\subset E^3\subset \cdots$ of compact subsets of $E$, having finite $\nu$-measure, such that $\bigcup_{l=1}^\infty E^l=E$. For every $l\in \mathbb{N}$, let $D^l:=\left\lbrace E^l_1,E^l_2,\ldots, E^l_{r_l}\right\rbrace\subset \mathcal{E}$ be a partition of $E^l$ finer than $\left\lbrace E^1,E^2\setminus E^1, \ldots , E^l\setminus E^{l-1}\right\rbrace$ such that for every $1\leq j \leq r_l$, the diameter of $E^l_j$ is less than $\varepsilon_l$. We suppose that $\varepsilon_l\downarrow 0$ as $l\to \infty$ for convergence of the numerical schemes. 

Concerning time discretization, we divide the time interval $[0,T]$ into $m$ subintervals of equal length and set $\delta_m=T/m, t_i=i\delta_m, i\in\{0,1,2,\ldots\}$. 
Now we wish to use these discretizations of time and the spaces $U, H, V, E$ to present explicit and implicit numerical schemes for equation \eqref{equation1}. 

\subsection{The Explicit Numerical Scheme} 

Let us first formulate $\tilde{A}^m$, $\tilde{B}^m$ and $\tilde{F}^{m,l}$ as approximations of operators $A$, $B$ and $F$. For all $x\in V$, all $\xi \in E$ and all $t\in (t_{i-1},t_i]$ we set 
\[\tilde{A}_t^m(x):=\tilde{A}_{t_i}^m(x),\quad\tilde{B}_t^m(x):=\tilde{B}_{t_i}^m(x),\quad\tilde{F}_t^{m,l}(x,\xi):=\tilde{F}_{t_i}^{m,l}(x,\xi),\] 
where 
\[\tilde{A}_{t_0}^m(x):=\tilde{A}_{t_1}^m(x):=0,\ \tilde{B}_{t_0}^m(x):=\tilde{B}_{t_1}^m(x):=0,\ \tilde{F}_{t_0}^{m,l}(x,\xi):=\tilde{F}_{t_1}^{m,l}(x,\xi):=0,\] 
\[\tilde{A}_{t_i}^m(x):=\frac{1}{\delta_m}\int_{t_{i-2}}^{t_{i-1}} A_s(x) \mathrm{d} s,\quad\tilde{B}_{t_i}^m(x):=\frac{1}{\delta_m}\int_{t_{i-2}}^{t_{i-1}} B_s(x) \mathrm{d} s,\quad 2\leq i\leq m.\] 
If $\xi\in E\setminus E^l$, we set for all $x\in V$ and $t \in [0,T]$, 
\[\tilde{F}_t^{m,l}(x,\xi):=0.\] 
If for some $1\leq j\leq r_l$, $\xi \in E^l_j$, we set 
\[\tilde{F}_{t_i}^{m,l}(x,\xi):=\frac{1}{\delta_m\cdot\nu(E^l_j)}\int_{t_{i-2}}^{t_{i-1}}\int_{E^l_j} F_s(x,\eta) \nu(\mathrm{d} \eta) \mathrm{d} s,\quad 2\leq i\leq m.\] 
The explicit discretization scheme is as follows 
\[u_{m,l}^n(t_0):=0, \quad u_{m,l}^n(t_1):=\Pi_n \zeta, \] 
\begin{equation}\label{explicit}
	\begin{split} 
		u^n_{m,l}(t_{i})&:=u^n_{m,l}(t_{i-1})+\delta_m\Pi_n\tilde{A}_{t_i}^m\left(u^n_{m,l}(t_{i-1})\right)\\&\quad+\Pi_n\tilde{B}_{t_i}^m\left(u^n_{m,l}(t_{i-1})\right)\left(W^l_{t_i}-W^l_{t_{i-1}}\right)\\ 
		&\quad+\int_E\Pi_n\tilde{F}_{t_i}^{m,l}\left(u^n_{m,l}(t_{i-1}),\xi\right)\tilde{N}\left((t_{i-1},t_i],\mathrm{d}\xi\right),\quad 2\leq i\leq m+1, 
	\end{split}
\end{equation} 
\[t_{i-1}<t\leq t_i: \quad u_{m,l}^n(t):=u^n_{m,l}(t_{i}),\quad 1\leq i\leq m.\] 
By using induction on $i$, it is clear that $u^n_{m,l}(t_{i})$ is $\mathcal{F}_{t_i}$-measurable. Every trajectory of $u_{m,l}^n$ is a left continuous step function. When $n, l$ and $m$ tend to infinity, the stochastic processes $u_{m,l}^n$ may not converge. Let us first introduce the following notation. 
\begin{definition} \label{normineq}
	For the space $V_n$, define 
		\[\mathcal{C}(n):=\mathrm{ess } \sup_{\substack{{t\in[0,T]}\\ {\omega\in \Omega}}}\sup_{v\in V_n}\frac{\left\lVert \Pi_n A_t(v)\right\rVert_H^2}{\left\lVert A_t(v)\right\rVert_{V^\star}^2}\] 
\end{definition} 

Now the following convergence theorem for the explicit scheme, with weaker condition than \cite[Theorem 2.8]{gyongy2005discretization}, holds. 
\begin{theorem} \label{explicitconvergence}
	Suppose conditions \ref{C1}-\ref{C7} with $p=2$. If $n, l$ and $m$ tend to infinity such that 
	\begin{equation}\label{condition}
	\mathcal{C}(n)\max_{1\leq i\leq m} \int_{t_{i-1}}^{t_i}\lambda(s)\mathrm{d} s\rightarrow 0,
	\end{equation}
	then the sequence of stochastic processes $u_{m,l}^n$ converges weakly in $\mathscr{L}_V^p(\lambda)$ to $u$, a solution of equation \eqref{equation1}. In addition $u_{m,l}^n(T)$ converges strongly in $L^2(\Omega;H)$ to $u(T)$. 
\end{theorem} 
In \cite[Theorem 2.8]{gyongy2005discretization}, $\sum_{k=1}^{l_n}\left\lVert e_k\right\rVert_V^2$ is considered instead of $\mathcal{C}(n)$ which has bigger order (see section \ref{example}). For $D = (0, 1)$, $V = W^{1,2}_0(D), H=L^2(D), 
	Au=\frac{\partial^2 u}{\partial x^2}$, and $e_n:=\sin (n\pi\cdot), n\in\mathbb{N}$ 
	the condition in \cite[Theorem 2.8]{gyongy2005discretization} reads as $\frac{n^3}{m}\to 0$ and \eqref{condition} reads as $n^2/m\to 0$. Note that by uniform continuity of the function $[0,T]\ni t\mapsto \int_0^t \lambda(s)\mathrm{d} s$, we have 
	\[ \max_{1\leq i\leq m} \int_{t_{i-1}}^{t_i}\lambda(s) \mathrm{d} s\to 0 \qquad\text{ as } m\to \infty.\]

\subsection{The Implicit Numerical Schemes} 

Here we discretize the operators $B$ and $F$ in the same way as in the explicit scheme. 
But for the operator $A$ we set the value of its discrete approximation $A^m$ at time $t$ to the average of $A$ over the subinterval containing $t$, instead of its preceding subinterval. More precisely, 
\[A_{t_0}^m(x):=0,\] 
\[A_{t_i}^m(x):=\frac{1}{\delta_m}\int_{t_{i-1}}^{t_i}A_s(x)\mathrm{d} s,\quad 1\leq i\leq m,\] 
\[t_{i-1}<t\leq t_i: \quad A_t^m:=A_{t_i}^m,\quad 1\leq i\leq m.\] 
With respect to the above introduced discretization of time, space $U$ and the measure $\nu$ we then define the following scheme 
\[u^{m,l}(t_0):=\zeta,\] 
\begin{equation}
	\label{implicitnotspace} 
	\begin{split}
		u^{m,l}(t_{i})&:=u^{m,l}(t_{i-1})+\delta_m A_{t_i}^m\left(u^{m,l}(t_{i})\right)+\tilde{B}_{t_i}^m\left(u^{m,l}(t_{i-1})\right)\left(W^l_{t_i}-W^l_{t_{i-1}}\right)\\&\quad+\int_E\tilde{F}_{t_i}^{m,l}\left(u^{m,l}(t_{i-1}),\xi\right)\tilde{N}\left((t_{i-1},t_i],\mathrm{d}\xi\right),\quad 1\leq i\leq m+1, 
	\end{split} 
\end{equation}
\[t_{i-1}<t\leq t_i:\quad u^{m,l}(t):=u^{m,l}(t_{i}),\quad 1\leq i\leq m.\] 
Adding the projection $\Pi_n$, we get another implicit scheme:
\[u^{n,m,l}(t_0):=\Pi_n\zeta, \] 
\begin{equation}\label{implicit} 
	\begin{split}
		u^{n,m,l}(t_{i})&:=u^{n,m,l}(t_{i-1})+\delta_m\Pi_n A_{t_i}^m\left(u^{n,m,l}
		(t_{i})\right)\\&\quad+\Pi_n\tilde{B}_{t_i}^m\left(u^{n,m,l}(t_{i-1})\right)\left(W^l_{t_i}-
		W^l_{t_{i-1}}\right)\\&\quad+\int_E\Pi_n\tilde{F}_{t_i}^{m,l}\left(u^{n,m,l}(t_{i-1}),
		\xi\right)\tilde{N}\left((t_{i-1},t_i],\mathrm{d}\xi\right),\quad 1\leq i\leq m+1, 
	\end{split} 
\end{equation}
\[t_{i-1}<t\leq t_i:\quad u^{n,m,l}(t):=u^{n,m,l}(t_{i}),\quad 1\leq i\leq m.\] 
Equations \eqref{implicitnotspace} and \eqref{implicit} have unique solutions 
$u^{m,l}(t_{i})$ and $u^{n,m,l}(t_{i})$ respectively, for $m$ sufficiently large.
This fact is stated in the next theorem which is similar to 
\cite[Theorem 2.9]{gyongy2005discretization}. 
\begin{theorem}\label{implicitexistance} 
	Assume conditions \ref{C1}-\ref{C7} with $p\in[2,\infty)$. Then there is a natural number $m_0\geq 1$ such that for every $m\geq m_0$ and $l\geq 1$, equation \eqref{implicitnotspace} has a unique solution $u^{m,l}(t_{i})$ that is $\mathcal{F}_{t_i}$-measurable and $\mathbb{E}\left\lVert u^{m,l}(t_{i})\right\rVert_V^p<\infty$ for each $i=1,2,\ldots,m$. Similarly there exists a natural number $m_0$ such that for every $m\geq m_0$ and $n,l\geq 1$, equation \eqref{implicit} has a unique solution $u^{n,m,l}(t_{i})$ that is $\mathcal{F}_{t_i}$-measurable and $\mathbb{E}\left\lVert u^{n,m,l}(t_{i})\right\rVert_V^p<\infty$ for each $i=1,2,\ldots,m$. 
\end{theorem} 
The convergence theorem for the implicit schemes, which is analogous to 
\cite[Theorem 2.10]{gyongy2005discretization}, is given as follows:

\begin{theorem}\label{implicitconvergence} 
	Assume \ref{C1}-\ref{C7} with $p\geq 2$. If $m$ and $l$ converge to infinity, then $u^{m,l}$ converges weakly in $\mathscr{L}^p_V(\lambda)$ to $u$, the solution of equation\eqref{equation1} and $u^{m,l}(T)$ converges strongly in ${L}^2(\Omega,H)$ to $u(T)$. Similarly, if $m,l$, and $n$ tend to infinity, $u^{n,m,l}$ converges weakly in $\mathscr{L}^p_V(\lambda)$ to $u$ and $u^{n,m,l}(T)$ converges strongly in ${L}^2(\Omega,H)$ to $u(T)$. 
\end{theorem}

\section{Proof of Results}

\subsection{Convergence of the Explicit Scheme}

First we obtain the integral form of equation \eqref{explicit}. For $t_{i}<t\leq t_{i+1},0\leq i\leq m-1$ set 
\[\kappa_1(t):=t_{i}\quad\kappa_2(t):=t_{i+1},\] 
and for $t_0$, $\kappa_1(t_0)=\kappa_2(t_0)=t_0$. Then 
\begin{equation}\label{ex int form}
	\begin{split} 
		u_{m,l}^n(t)&=\Pi_n\zeta \mathbf{1}_{\left\lbrace t>t_0\right\rbrace}+\int_0^{\kappa_1(t)}\Pi_nA_s\left(u_{m,l}^n(s)\right)\mathrm{d} s\\&\quad+\int_0^{\kappa_2(t)}\Pi_n\tilde{B}^m_s\left(u_{m,l}^n\left(\kappa_1(s)\right)\right)\tilde{\Pi}_l\mathrm{d} W_s\\ 
		& \quad+\int_0^{\kappa_2(t)}\int_E\Pi_n\tilde{F}^{m,l}_s\left(u_{m,l}^n\left(\kappa_1(s)\right)\right)\tilde{N}(\mathrm{d} s,\mathrm{d}\xi). 
	\end{split}
\end{equation} 
We wish to prove boundedness of $u_{m,l}^n$ and the integrands of the above equation in some convenient reflexive Banach spaces. Then weakly compactness of bounded sequences in reflexive Banach spaces (see e.g. \cite[Theorem 3.18]{brezis2010functional}) implies weak convergence of some subsequences of $u_{m,l}^n$ and the integrands to some stochastic processes, like $\bar{u}_\infty, a_{\infty}, b_{\infty}$ and $f_{\infty}$, where for all $t \in [0,T]$ and all $z \in V$, 
\begin{align*}(\bar{u}_\infty(t),z)_H&=(\zeta ,z)_H+\int_0^t \langle a_{\infty}(s),z\rangle \mathrm{d} s+\int_0^t \left(b_{\infty}(s)\mathrm{d} W_s,z\right)_H\\&\quad+\int_0^t\int_E(f_{\infty}(s,\xi),z)_H\tilde{N}(\mathrm{d} s,\mathrm{d}\xi) \quad a.s. 
\end{align*} 
We will get that $(\zeta ,a_{\infty}, b_{\infty},f_{\infty})$ belongs to the set $\mathcal{A}$ and $I_y(\zeta ,a_{\infty}, b_{\infty},f_{\infty})\leq0$ for all $y\in\mathscr{G}_\mathcal{BF}$. So by Theorem \ref{characterize}, $\bar{u}_\infty$ will have a modification which is a solution of equation \eqref{equation1}. This will complete the proof. 

%

\begin{proof}[Proof of Theorem \ref{explicitconvergence}]
	The first step asserts the boundedness of $u_{m,l}^n$ and also the integrands of equation \eqref{ex int form}. 
	\paragraph*{\textbf{Step 1. }}\label{boundedness} 
	\textit{
		If $0<\gamma<1$ and 
		\[I_\gamma:=\left\lbrace (n,m,l):\alpha\, C_\mathpzc{B}(n)\,\max_{1\leq i\leq m}\int_{t_{i-1}}^{t_i}\lambda(s)\mathrm{d} s\leq \gamma\right\rbrace,\] 
		then the following functions of $(n, m, l)$ are bounded on $I_\gamma$:
		\begin{enumerate}[label=(\roman{*})] 
			\item \label{boundedness:i}$\sup_{t\in [0,T]}\mathbb{E}\left\lVert u_{m,l}^n(t)\right\rVert_H^2 $, 
			\item \label{boundedness:ii}$\mathbb{E}\int_0^T\left\lVert u_{m,l}^n(t)\right\rVert_V^2\lambda(t)\mathrm{d} t$, 
			\item \label{boundedness:iii}$\mathbb{E}\int_0^T\left\lVert A_t\left(u_{m,l}^n(t)\right)\right\rVert_{V^*}^2\lambda(t)^{-1}\mathrm{d} t$,
			\item \label{boundedness:iv}$\mathbb{E}\int_0^T\left\lVert\Pi_n\tilde{B}^m_t\left(u_{m,l}^n(\kappa_1(t))\right)\tilde{\Pi}_l\right\rVert_2^2\mathrm{d} t$, 
			\item \label{boundedness:v}$\mathbb{E}\int_0^T\int_E\left\lVert\Pi_n\tilde{F}^{m,l}_t\left(u_{m,l}^n(\kappa_1(t)),\xi\right)\right\rVert_{H}^2\nu(\mathrm{d}\xi)\mathrm{d} t$. 
	\end{enumerate} }
	\paragraph*{\text{Proof of Step 1. }} 
	By the definition of the explicit scheme, i.e. equation \eqref{explicit}, for \linebreak$1\leq i\leq m$, we have 
	\begin{equation}\label{before15}
		\begin{split} 
			\mathbb{E}\left\lVert u^n_{m,l}(t_{i})\right\rVert_H^2 &= \mathbb{E}\left\lVert u^n_{m,l}(t_{i-1})\right\rVert_H^2 +\delta_m^2\mathbb{E}\left\lVert\Pi_n\tilde{A}^m_{t_{i}}(u_{m,l}^n(t_{i-1}))\right\rVert_H^2 \\ 
			&\quad+\delta_m\mathbb{E}\left\lVert\Pi_n\tilde{B}^m_{t_i}\left(u_{m,l}^n(t_{i-1})\right)\tilde{\Pi}_l\right\rVert_2^2 \\&\quad+\delta_m\mathbb{E}\int_E\left\lVert\Pi_n\tilde{F}^{m,l}_{t_i}(u_{m,l}^n(t_{i-1}),\xi)\right\rVert_H^2\nu(\mathrm{d}\xi)\\&\quad+2\delta_m\mathbb{E}\left\langle u_{m,l}^n(t_{i-1}),\Pi_n \tilde{A}_{t_i}^m(u_{m,l}^n(t_{i-1}))\right\rangle. \end{split} 
	\end{equation}
	Note that, since $W_{t_i}-W_{t_{i-1}}$ and $\tilde{N}((t_{i-1},t_i],\mathrm{d} \xi)$ are independent of each other and independent of $\mathcal{F}_{t_{i-1}}$, the expectation of their cross product with any $\mathcal{F}_{t_{i-1}}$-adapted random variable is zero. By using the definition of the discretized operators $\tilde{A}^m,\tilde{B}^m$ and $\tilde{F}^{m,l}$, we get 
	\begin{align*} 
		\mathbb{E}\left\lVert u^n_{m,l}(t_{i})\right\rVert_H^2 &\leq \mathbb{E}\left\lVert u^n_{m,l}(t_{i-1})\right\rVert_H^2+\mathbb{E}\left\lVert\int_{t_{i-2}}^{t_{i-1}}\Pi_nA_s(u^n_{m,l}(t_{i-1}))\mathrm{d} s\right\rVert_{H}^2\\& +\delta_m^{-1}\mathbb{E}\left\lVert\int_{t_{i-2}}^{t_{i-1}}\Pi_nB_s(u^n_{m,l}(t_{i-1}))\tilde{\Pi}_l\mathrm{d} s\right\rVert_2^2\\ 
		& + \delta_m^{-1}\sum_{1\leq j\leq r_l}\nu(E^l_j)^{-1}\mathbb{E}\left\lVert\int_{t_{i-2}}^{t_{i-1}}\int_{E^l_j}\Pi_nF_s(u^n_{m,l}(t_{i-1}),\xi)\nu(\mathrm{d} \xi)\mathrm{d} s\right\rVert_H^2\\&+2\mathbb{E}\left\langle u^n_{m,l}(t_{i-1}),\int_{t_{i-2}}^{t_{i-1}}A_s(u^n_{m,l}(t_{i-1}))\mathrm{d} s\right\rangle \end{align*} 
	and by using Definition \ref{normineq} for $\left\lVert\Pi_nA_s(u^n_{m,l}(t_{i-1}))\right\rVert_H^2$ and the Cauchy-Schwartz's inequality 
	\begin{equation}\label{ghablan(4)}
		\begin{split}
			\MoveEqLeft[1]\mathbb{E}\left\lVert u^n_{m,l}(t_{i})\right\rVert_H^2\leq \mathbb{E}\left\lVert u^n_{m,l}(t_{i-1})\right\rVert_H^2\\ & +\left(\max_{1\leq i\leq m}\int_{t_{i-1}}^{t_i}\lambda(s)\mathrm{d} s\right)\mathcal{C}(n)\mathbb{E}\int_{t_{i-2}}^{t_{i-1}}\left\lVert A_s(u^n_{m,l}(t_{i-1}))\right\rVert_{V^*}^2\lambda^{-1}(s)\mathrm{d} s\\&+\mathbb{E}\int_{t_{i-2}}^{t_{i-1}}\Bigg[\left\lVert\Pi_nB_s(u^n_{m,l}(t_{i-1}))\tilde{\Pi}_l\right\rVert_2^2 + \int_E\left\lVert \Pi_nF_s(u^n_{m,l}(t_{i-1}),\xi)\right\rVert_H^2\nu(\mathrm{d}\xi)\\ &+2\left\langle u^n_{m,l}(t_{i-1}),A_s(u^n_{m,l}(t_{i-1}))\right\rangle\Bigg] \mathrm{d} s. 
		\end{split} 
	\end{equation}
	The coercivity condition \ref{C2} and the growth condition \ref{C3} yield that the right hand side of \eqref{ghablan(4)} is less than or equal to
	\begin{align*} 
		&\mathbb{E}\left\lVert u^n_{m,l}(t_{i-1})\right\rVert_H^2 +\left(\max_{1\leq i\leq m}\int_{t_{i-1}}^{t_i}\lambda(s)\mathrm{d} s\right)\mathcal{C}(n)\\&\times\mathbb{E}\int_{t_{i-2}}^{t_{i-1}}\left[\alpha\lambda(s)\left\lVert u^n_{m,l}(t_{i-1})\right\rVert_{V}^2+ K_2(s)\right]\mathrm{d} s \\ 
		&+\mathbb{E}\int_{t_{i-2}}^{t_{i-1}}\left[-\lambda(s)\left\lVert u^n_{m,l}(t_{i-1})\right\rVert_{V}^2+ \mu(s)\left\lVert u^n_{m,l}(t_{i-1})\right\rVert_{H}^2+K_1(s)\right] \mathrm{d} s.
	\end{align*} 
	Now define $\rho:=1-\alpha\left(\max_{1\leq i\leq m}\int_{t_{i-1}}^{t_i}\lambda(s)\mathrm{d} s\right)\mathcal{C}(n)$. Since $(n,m,l)\in I_\gamma$, we have $\rho>0$ and 
	\begin{align*} 
		\MoveEqLeft[3]\mathbb{E}\left\lVert u^n_{m,l}(t_{i})\right\rVert_H^2 +\rho\ \mathbb{E}\int_{t_{i-2}}^{t_{i-1}}\lambda(s)\left\lVert u^n_{m,l}(t_{i-1})\right\rVert_{V}^2 \mathrm{d} s \\ &\leq \mathbb{E}\left\lVert u^n_{m,l}(t_{i-1})\right\rVert_H^2 +\mathbb{E}\left\lVert u^n_{m,l}(t_{i-1})\right\rVert_H^2 \int_{t_{i-2}}^{t_{i-1}}\mu(s)\mathrm{d} s\\&\quad+\mathbb{E}\int_{t_{i-2}}^{t_{i-1}}\left[K_1(s)+K_2(s)/\alpha\right] \mathrm{d} s. 
	\end{align*} 
	Summing up the above inequality for $i=1,2,\ldots,k$, with $1\leq k\leq m+1$, we get 
	\begin{equation}\label{boundfortwo} 
		\mathbb{E}\left\lVert u_{m,l}^n(t_{k})\right\rVert_H^2 +\rho\ \mathbb{E}\int_{0}^{t_{k-1}}\lambda(s)\left\lVert u_{m,l}^n(s)\right\rVert_{V}^2 \mathrm{d} s \leq C+ \sum_{i=1}^{k-1}\alpha_i \mathbb{E}\left\lVert u^n_{m,l}(t_{i-1})\right\rVert_H^2, 
	\end{equation} 
	where $\alpha_i=\int_{t_{i-1}}^{t_{i}}\mu(s)\mathrm{d} s$ and $C=\mathbb{E}\left\lVert \zeta \right\rVert_H^2+\mathbb{E}\int_{0}^{T}\left[K_1(s)+K_2(s)/\alpha\right] \mathrm{d} s$. Now we neglect the second term on the left hand side of inequality above, and using induction and the fact that 
	\[\mathbb{E}\left\lVert u_{m,l}^n(t_{0})\right\rVert_H^2=0,\qquad\mathbb{E}\left\lVert u_{m,l}^n(t_{1})\right\rVert_H^2=\mathbb{E}\left\lVert\Pi_n\zeta \right\rVert_H^2\leq C,\] 
	we get the following inequality for $0\leq k \leq m$ 
	\begin{align*} 
		\mathbb{E}\left\lVert u_{m,l}^n(t_k)\right\rVert_H^2 &\leq C(1+\alpha_1)(1+\alpha_2)\cdots(1+\alpha_{k-1})\\ 
		&\leq C(1+\alpha_1)(1+\alpha_2)\cdots(1+\alpha_{m})\\ 
		&\leq C\left(1+\frac{\int_0^T \mu(s)\mathrm{d} s}{m}\right)^m \, . 
	\end{align*} 
	The sequence $C\left(1+\frac{\int_0^T \mu(s)\mathrm{d} s}{m}\right)^m, m\in\mathbb{N}$ converges to the finite number 
	\[
	C\exp\left( \int_0^T \mu(s)\mathrm{d} s\right), 
	\] 
	as $m\to\infty$, so we have
	\[
	\sup_{(n,m)\in I_\gamma}\sup_{t\in[0,T]}\mathbb{E}\left\lVert u_{m,l}^n(t)\right\rVert_H^2<\infty.
	\] 
	Equation \eqref{boundfortwo} for $k=m+1$ implies boundedness of \ref{boundedness:ii} over $I_\gamma$, which together with the growth condition \ref{C3} of $A$ implies 
	\begin{align*}
	\MoveEqLeft\sup_{(n,m)\in I_\gamma}\mathbb{E}\int_0^T\left\lVert A_t\left(u_{m,l}^n(t)\right)\right\rVert_{V^*}^2\lambda(t)^{-1}\mathrm{d} t\\&\leq\sup_{(n,m)\in I_\gamma}\mathbb{E}\int_0^T\left[\alpha\lambda(t)\left\lVert u_{m,l}^n(t)\right\rVert_{V}^2+K_2(t)\right]\mathrm{d} t<\infty,\end{align*} 
	so \ref{boundedness:iii} is also bounded over $I_\gamma$. By the definition of $\tilde{B}^m$ and Cauchy-Schwartz's inequality we have 
	\begin{align*} 
		\MoveEqLeft\mathbb{E}\int_0^T\left\lVert\Pi_n\tilde{B}^m_s\left(u_{m,l}^n(\kappa_1(s))\right)\tilde{\Pi}_l\right\rVert_2^2\mathrm{d} s \\&=\delta_m^{-1}\sum_{i=2}^m\mathbb{E}\left\lVert\int_{t_{i-2}}^{t_{i-1}}\Pi_nB_s(u^n_{m,l}(t_{i-1}))\tilde{\Pi}_l\mathrm{d} s\right\rVert_2^2\\&\leq\sum_{i=2}^m\mathbb{E}\int_{t_{i-2}}^{t_{i-1}}\lVert\Pi_nB_s(u^n_{m,l}(t_{i-1}))\tilde{\Pi}_l\rVert_2^2\mathrm{d} s\\&\leq\mathbb{E}\int_0^T\left\lVert\Pi_nB_s \left(u_{m,l}^n(s)\right)\tilde{\Pi}_l\right\rVert_2^2\mathrm{d} s. 
	\end{align*} 
	Similarly 
	\begin{align*}\MoveEqLeft\mathbb{E}\int_0^T\int_E\left\lVert\Pi_n\tilde{F}^{m,l}_t\left(u_{m,l}^n(\kappa_1(t)),\xi\right)\right\rVert_{H}^2\nu(\mathrm{d}\xi)\mathrm{d} t\\&\leq\mathbb{E}\int_0^T\int_E\left\lVert\Pi_nF_t\left(u_{m,l}^n(t),\xi\right)\right\rVert_{H}^2\nu(\mathrm{d}\xi)\mathrm{d} t. 
	\end{align*} 
	By using Proposition \ref{proposition1}, we get 
	\begin{align*}&\sup_{(n,m,l)\in I_\gamma}\mathbb{E}\int_0^T\left[\left\lVert\Pi_nB_s \left(u_{m,l}^n(s)\right)\tilde{\Pi}_l\right\rVert_2^2 +\int_E\left\lVert\Pi_nF_s\left(u_{m,l}^n(s),\xi\right)\right\rVert_{H}^2\nu(\mathrm{d}\xi)\right]\mathrm{d} s\\&\leq\sup_{(n,m,l)\in I_\gamma}\mathbb{E}\int_0^T\left[6\left(\alpha^{1/q}+1/p\right)\lambda(s)\left\lVert u_{m,l}^n(s)\right\rVert_V^2+\frac{16}{q}K_2(s)+2K_1(s)\right]\mathrm{d} s\\&<\infty\, . 
	\end{align*} 
	Hence, \ref{boundedness:iv} and \ref{boundedness:v} are bounded too, and the proof of Step 1 is completed. 
	
	\paragraph{\textbf{Step 2. }}
	\textit{
		Let $(n,m,l)$ be a sequence from $I_\gamma$ for some $\gamma\in(0,1)$, such that $m$, $n$ and $l$ converge to infinity. Then it contains a subsequence, denoted also by $(n,m,l)$, such that 
		\begin{enumerate}[label=(\roman{*})] 
			\item \label{converg of explicit:i}$u_{m,l}^n$ converges weakly in $\mathscr{L}_V^2(\lambda)$ to some progressively measurable process $\bar{u}_\infty$, 
			\item \label{converg of explicit:ii}$u_{m,l}^n(T)$ converges weakly in $L^2(\Omega;H)$ to some random variable $u_{T\infty}$, 
			\item \label{converg of explicit:iii}$A_\cdot \left(u_{m,l}^n(\cdot)\right)$ converges weakly in $\mathscr{L}_{V^*}^2(\lambda^{-1})$ to some progressively measurable process $a_\infty$, 
			\item \label{converg of explicit:iv}$\Pi_n\tilde{B}_\cdot^m\left(u_{m,l}^n(\kappa_1(\cdot))\right)\tilde{\Pi}_l$ converges weakly in $\mathscr{L}_{L_2(U,H)}^2(\mathcal{BF})$ to some process $b_\infty$, 
			\item \label{converg of explicit:v}$\Pi_n\tilde{F}_\cdot^{m,l}\left(u_{m,l}^n(\kappa_1(\cdot)),\star\right)$ converges weakly in $L^2([0,T]\times\Omega\times E, \mathcal{P}\otimes\mathcal{E},\mathrm{d} t\otimes\mathbb{P}\otimes\nu;H)$ to some process $f_\infty$, 
			\item \label{converg of explicit:vi}$(\zeta ,a_\infty,b_\infty,f_\infty)\in \mathcal{A}$ and for all $z\in V$, and $\mathrm{d} t\otimes \mathbb{P}$- almost all $(t,\omega)$ we have 
			\begin{equation}\label{v infty} 
				\begin{split}
					(\bar{u}_\infty(t),z)_H&=(\zeta ,z)_H+\int_0^t \langle a_{\infty}(s),z\rangle \mathrm{d} s+\int_0^t \left(b_{\infty}(s)\mathrm{d} W_s,z\right)_H\\&\quad+\int_0^t\int_E(f_{\infty}(s,\xi),z)_H\tilde{N}(\mathrm{d} s,\mathrm{d}\xi) 
				\end{split} 
			\end{equation}
			and for all $z\in V$, almost surely 
			\begin{equation}\label{u infty} 
				\begin{split}
					(u_{T\infty},z)_H&=(\zeta ,z)_H+\int_0^T\langle a_\infty(s),z\rangle \mathrm{d} s+\int_0^T\left(z,b_\infty(s)\mathrm{d} W_s\right)_H\\&\quad+\int_0^T\int_E(f_\infty(s,\xi),z)_H\tilde{N}(\mathrm{d} s,\mathrm{d}\xi) \, . 
				\end{split}
			\end{equation} 
	\end{enumerate} }
	
	\paragraph{\text{Proof of Step 2. }}
	The convergences in \ref{converg of explicit:i}-\ref{converg of explicit:v} can be immediately concluded from Step 1, except the fact that $\bar{u}_\infty$ and $a_\infty$ are progressively measurable. Note that $\bar{u}_\infty(t)$ and $a_\infty(t)$ are $\mathcal{F}_{t+\delta_m}$-adapted processes for each $m\geq 1$, so they are $\mathcal{F}_t$-adapted and also $\mathcal{B}([0,T])\otimes\mathcal{F}$-measurable. Hence they have progressively measurable modifications, that will replace them in the following (see e.g. \cite{ondrejat2013existence}). It remains to prove \ref{converg of explicit:vi}. Fix $N\in\mathbb{N}$. It is sufficient to verify \ref{converg of explicit:vi} for $z\in V_N$ because $\bigcup_{N=1}^\infty V_N$ is dense in $V$. Both sides of \eqref{v infty} belong to the Hilbert space $\mathscr{L}_\mathbb{R}^2$. Therefore to verify \eqref{v infty}, it is sufficient to prove that the inner products of both sides and any $\varphi\in \mathscr{L}_\mathbb{R}^2$ are the same, i.e.
	\begin{equation}\label{effect of v infty} 
		\begin{split}
			\mathbb{E}\int_0^T \left(\bar{u}_\infty(t),z\right)_H\varphi(t)\mathrm{d} t&=\mathbb{E}\int_0^T\left(\zeta ,z\right)_H\varphi(t)\mathrm{d} t\\&\quad+\mathbb{E}\int_0^T\left(\int_0^t\left\langle a_\infty(s),z\right\rangle \mathrm{d} s\right)\varphi(t)\mathrm{d} t\\& \quad+\mathbb{E}\int_0^T\left(\int_0^t\left(b_\infty(s)\mathrm{d} W_s,z\right)_H\right)\varphi(t)\mathrm{d} t\\&\quad +\mathbb{E}\int_0^T\left(\int_0^t\int_E \left(f_\infty(s,\xi),z\right)_H\tilde{N}(\mathrm{d} s,\mathrm{d}\xi)\right)\varphi(t)\mathrm{d} t \, . 
		\end{split}
	\end{equation} 
	Note that the integral form of the explicit scheme \eqref{ex int form} yields for $z\in V_N$ and $n\geq N$ 
	\begin{align*} 
		(u_{m,l}^n(t),z)_H&=(\zeta ,z)_H\mathbf{1}_{\left\lbrace t>t_0\right\rbrace}+\int_0^{\kappa_1(t)}\left\langle A_s\left(u_{m,l}^n(s)\right),z\right\rangle\mathrm{d} s\\&\quad+\int_0^{\kappa_2(t)}\left(z,\Pi_n\tilde{B}^m_s\left(u_{m,l}^n\left(\kappa_1(s)\right)\right)\mathrm{d} W^l_s\right)_H\\ 
		&\quad+\int_0^{\kappa_2(t)}\int_E\left(\Pi_n\tilde{F}^{m,l}_s\left(u_{m,l}^n\left(\kappa_1(s)\right)\right),z\right)_H\tilde{N}(\mathrm{d} s,\mathrm{d}\xi)\quad a.s. 
	\end{align*} 
	Taking the inner products of both sides and $\varphi$, we get for $n\geq N$, 
	\begin{equation}\label{eq:Ji Ri's} 
		\begin{split}
			\mathbb{E}\int_0^T\left( u_{m,l}^n(t),z\right)_H\varphi(t) \mathrm{d} t&=\mathbb{E}\int_0^T\left(\zeta ,z\right)_H\varphi(t)\mathrm{d} t+J_1+J_2+J_3\\&\quad-R_1-R_2-R_3 \, , 
		\end{split} 
	\end{equation}
	where 
	\begin{align*} 
		J_1&=\mathbb{E}\int_0^T\varphi(t)\left(\int_0^t \left\langle A_s\left(u_{m,l}^n(s)\right),z\right\rangle \mathrm{d} s\right)\mathrm{d} t,\\ 
		J_2&=\mathbb{E}\int_0^T\varphi(t)\int_0^t\left( \Pi_n\tilde{B}^m_s\left(u_{m,l}^n(\kappa_1(s))\right)\tilde{\Pi}_l\mathrm{d} W_s,z\right)_H\mathrm{d} t,\\ 
		J_3&=\mathbb{E}\int_0^T\varphi(t)\left(\int_0^t\int_E \left(\Pi_n\tilde{F}^{m,l}_s\left(u_{m,l}^n(\kappa_1(s)),\xi\right),z\right)_H\tilde{N}(\mathrm{d} s,\mathrm{d}\xi)\right)\mathrm{d} t,\end{align*}
	and
	\begin{align*}
		R_1&=\mathbb{E}\int_0^T\varphi(t)\left(\int_{\kappa_1(t)}^t \left\langle A_s\left(u_{m,l}^n(s)\right),z\right\rangle \mathrm{d} s\right)\mathrm{d} t,\\ 
		R_2&=\mathbb{E}\int_0^T\varphi(t)\int_{t}^{\kappa_2(t)} \left(\Pi_n\tilde{B}^m_s\left(u_{m,l}^n(\kappa_1(s))\right)\tilde{\Pi}_l\mathrm{d} W_s,z\right)_H\mathrm{d} t,\\ 
		R_3&=\mathbb{E}\int_0^T\varphi(t)\left(\int_{t}^{\kappa_2(t)}\int_E \left(\Pi_n\tilde{F}^{m,l}_s\left(u_{m,l}^n(\kappa_1(s)),\xi\right),z\right)_H\tilde{N}(\mathrm{d} s,\mathrm{d}\xi)\right)\mathrm{d} t.\\ 
	\end{align*} 
	Our goal is to identify the limits of $J_i$'s and $R_i$'s. For $J_1$, consider the linear operator $S_1:\mathscr{L}_{V^*}^q(\lambda^{1-q})\to \mathscr{L}_\mathbb{R}^q$, defined by 
	\[S_1(g)(t):=\int_0^t\left\langle g(s),z\right\rangle \mathrm{d} s \] 
	for all $g\in \mathscr{L}_{V^*}^q(\lambda^{1-q})$. $S_1$ is bounded, because by H\"older's inequality we have that 
	\begin{align*} 
		\MoveEqLeft[3]\mathbb{E}\int_0^T\left\lvert S_1(g)(t)\right\rvert^q \mathrm{d} t \\&=\mathbb{E}\int_0^T\left\lvert\int_0^t\left\langle g(s),z\right\rangle \mathrm{d} s\right\rvert^q\mathrm{d} t\\&\leq \mathbb{E}\int_0^T\left(\int_0^t \left\rvert\left\langle g(s),z\right\rangle \right\rvert^q\lambda(s)^{-q/p} \mathrm{d} s\right)\left(\int_0^t\lambda(s) \mathrm{d} s\right)^{q/p}\mathrm{d} t\\&\leq \|z\|^q_V\mathbb{E}\int_0^T\left(\int_0^T \left\lVert g(s)\right\rVert _{V^*}^q\lambda(s)^{1-q} \mathrm{d} s\right)\left(\int_0^T\lambda(s) \mathrm{d} s\right)^{q/p}\mathrm{d} t \\ &\leq \|z\|^q_V T \left(\int_0^T\lambda(s) \mathrm{d} s\right)^{q/p} \left\lvert g\right\rvert_{\mathscr{L}_{V^*}^q(\lambda^{1-q},\mathcal{BF})}^q \, . 
	\end{align*} 
	So, $S_1$ is continuous with respect to the weak topologies. Thus by \ref{converg of explicit:iii}, i.e. ,
	\[A_\cdot\left(u_{m,l}^n(\cdot)\right)\rightharpoonup a_\infty \quad \text{in } \mathscr{L}_{V^*}^q(\lambda^{1-q}),\] 
	we obtain that 
	\[S_1\left(A_\cdot\left(u_{m,l}^n(\cdot)\right)\right)\rightharpoonup S_1\left( a_\infty\right) \quad \text{in } \mathscr{L}_\mathbb{R}^q,\] 
	therefore 
	\begin{align*} 
		J_1&=\mathbb{E}\int_0^T\varphi(t)\left(\int_0^t \left\langle A_s\left(u_{m,l}^n(s)\right),z\right\rangle \mathrm{d} s\right)\mathrm{d} t\\&\to \mathbb{E}\int_0^T\varphi(t)\left(\int_0^t \left\langle a_\infty(s),z\right\rangle \mathrm{d} s\right)\mathrm{d} t. 
	\end{align*}
	Now for $J_2$, take $S_2$, the bounded linear operator as follows: 
	\[
	S_2:\mathscr{L}_{L_2(U,H)}^2(\mathcal{BF})\to \mathscr{L}_\mathbb{R}^2 
	\] 
	\[
	S_2(g)(t)=\int_0^t \left(z,g(s)\mathrm{d} W_s\right)_H \, . 
	\] 
	The boundedness of $S_2$ yields that $S_2$ is continuous with respect to the weak topologies. Therefore, by using \[\Pi_n\tilde{B}_\cdot^m\left(u_{m,l}^n(\kappa_1(\cdot))\right)\tilde{\Pi}_l\rightharpoonup b_\infty \qquad \text{in}\ \mathscr{L}_{L_2(U,H)}^2(\mathcal{BF}),\]
	we obtain that 
	\begin{align*}
		J_2&=\mathbb{E}\int_0^T\varphi(t)\int_0^t\left( \Pi_n\tilde{B}^m_s\left(u_{m,l}^n(\kappa_1(s))\right)\tilde{\Pi}_l\mathrm{d} W_s,z\right)_H\mathrm{d} t\\&\to\mathbb{E}\int_0^T\varphi(t)\int_0^t\left( b_\infty(s)\mathrm{d} W_s,z\right)_H\mathrm{d} t. 
	\end{align*} 
	Similarly, let us define the linear operator $S_3$ as 
	\[S_3:L^2([0,T]\times\Omega\times E, \mathcal{P}\otimes\mathcal{E},\mathrm{d} t\otimes\mathbb{P}\otimes\nu;H)\to \mathscr{L}_\mathbb{R}^2\] 
	\[S_3(g)(t):=\int_0^t\int_E \left(g(s,\xi),z\right)_H\tilde{N}(\mathrm{d} s,\mathrm{d}\xi).\] 
	We have by It\^o-L{\'e}vy's isometry
	\begin{align*}
		\mathbb{E}\int_0^T \left\lVert S_3(g)(t)\right\rVert^2\mathrm{d} t&\leq \mathbb{E}\int_0^T \left\lVert \int_0^t\int_E \left(g(s,\xi),z\right)_H\tilde{N}(\mathrm{d} s,\mathrm{d}\xi)\right\rVert^2\mathrm{d} t\\&=\mathbb{E}\int_0^T \int_0^t\int_E \left(g(s,\xi),z\right)^2_H\nu(\mathrm{d}\xi) \mathrm{d} s\,\mathrm{d} t\\&\leq \left\lVert z\right\rVert_H^2\,T\mathbb{E}\int_0^T \int_E\left\lVert g(s,\xi)\right\rVert_H^2 \nu(\mathrm{d} \xi)\mathrm{d} s. 
	\end{align*}
	So $S_3$ is bounded linear operator and therefore it is continuous with respect to the weak topologies. Since
	\[\Pi_n\tilde{F}^{m,l}_\cdot\left(u_{m,l}^n(\kappa_1(\cdot)),*\right)\rightharpoonup f_\infty\, ,\]
	we get
	\begin{align*}
		J_3&=\mathbb{E}\int_0^T\varphi(t)\left(\int_0^t\int_E \left(\Pi_n\tilde{F}^{m,l}_s\left(u_{m,l}^n(\kappa_1(s)),\xi\right),z\right)_H\tilde{N}(\mathrm{d} s,\mathrm{d}\xi)\right)\mathrm{d} t\\&\to\mathbb{E}\int_0^T\varphi(t)\left(\int_0^t\int_E \left(f_\infty(s,\xi),z\right)_H\tilde{N}(\mathrm{d} s,\mathrm{d}\xi)\right)\mathrm{d} t.
	\end{align*} 
	Now we wish to prove that the "$R_i$"s tend to zero. Applying Cauchy-Schwartz's inequality yields 
	\begin{align*} 
		R_1^2 &=\left\lvert\mathbb{E}\int_0^T\varphi(t)\left(\int_{\kappa_1(t)}^t \left\langle A_s\left(u_{m,l}^n(s)\right),z\right\rangle \mathrm{d} s\right)\mathrm{d} t\right\rvert^2\\ 
		&\leq \left\lVert \varphi\right\rVert_{\mathscr{L}_\mathbb{R}^2}^2\mathbb{E}\int_0^T\left\lvert\int_{\kappa_1(t)}^t \left\langle A_s\left(u_{m,l}^n(s)\right),z\right\rangle \mathrm{d} s\right\rvert^2\mathrm{d} t
		\\&\leq \left(\max_{1\leq i\leq m}\int_{t_{i-1}}^{t_i}\lambda(s)\mathrm{d} s\right)\\&\quad \times\left\lVert z\right\rVert^2_V\left\lVert \varphi\right\rVert_{\mathscr{L}_\mathbb{R}^2}^2\mathbb{E}\int_0^T
			\int_{\kappa_1(t)}^t \left\lVert A_s\left(u_{m,l}^n(s)\right)\right\rVert_{V^*}^2\lambda^{-1}(s) \mathrm{d} s\,\mathrm{d} t\\&\leq \left(\max_{1\leq i\leq m}\int_{t_{i-1}}^{t_i}\lambda(s)\mathrm{d} s\right)\\&\quad\times \left\lVert z\right\rVert^2_V\left\lVert \varphi\right\rVert_{\mathscr{L}_\mathbb{R}^2}^2\mathbb{E}\int_0^T
			\int_{0}^T \left\lVert A_s\left(u_{m,l}^n(s)\right)\right\rVert_{V^*}^2\lambda^{-1}(s) \mathrm{d} s\,\mathrm{d} t
			\\&=\left(\max_{1\leq i\leq m}\int_{t_{i-1}}^{t_i}\lambda(s)\mathrm{d} s\right)\\&\quad\times \left\lVert z\right\rVert^2_V \ \left\lVert \varphi\right\rVert_{\mathscr{L}_\mathbb{R}^2}^2\ T \times\mathbb{E}\int_0^T \left\lVert A_s\left(u_{m,l}^n(s)\right)\right\rVert_{V^*}^2\lambda(s)^{-1}\mathrm{d} s.
	\end{align*} 
	Hence,  by (iii) of step 1, we get $R_1\to 0$ when $(n,m,l)\in I_\gamma$ for $0<\gamma<1$, $n\geq N$ and $m\to \infty$.
	For $R_2$, using Cauchy-Schwartz's inequality and It\^o's isometry, we have 
	\begin{align*} 
		R_2^2 & =\left\lvert \mathbb{E}\int_0^T\varphi(t)\int_{t}^{\kappa_2(t)} \left(\Pi_n\tilde{B}^m_s\left(u_{m,l}^n(\kappa_1(s))\right)\tilde{\Pi}_l\mathrm{d} W_s,z\right)_H\mathrm{d} t\right\rvert^2\\ 
		&\leq\left\lVert\varphi\right\rVert_{\mathscr{L}_\mathbb{R}^2}^2\mathbb{E}\int_0^T\left\lvert\int_{t}^{\kappa_2(t)}\left(\Pi_n\tilde{B}^m_s\left(u_{m,l}^n(\kappa_1(s))\right)\tilde{\Pi}_l\mathrm{d} W_s,z\right)_H\right\rvert^2\mathrm{d} t\\&\leq\left\lVert z\right\rVert_V^2\left\lVert \varphi\right\rVert_{\mathscr{L}_\mathbb{R}^2}^2\mathbb{E}\int_0^T\int_{t}^{\kappa_2(t)} \left\lVert \Pi_n\tilde{B}^m_s\left(u_{m,l}^n(\kappa_1(s))\right)\tilde{\Pi}_l\right\rVert_{2}^2 \mathrm{d} s\,\mathrm{d} t\\&=
		\left\lVert z\right\rVert^2_V \left\lVert \varphi\right\rVert_{\mathscr{L}_\mathbb{R}^2}^2\mathbb{E}\int_0^T \left(\kappa_2(t)-t\right)\left\lVert \Pi_n\tilde{B}^m_t\left(u_{m,l}^n(\kappa_1(t))\right)\tilde{\Pi}_l\right\rVert_{2}^2\mathrm{d} t\\&\leq \delta_m 
		\left\lVert z\right\rVert^2_V \left\lVert \varphi\right\rVert_{\mathscr{L}_\mathbb{R}^2}^2\mathbb{E}\int_0^T \left\lVert \Pi_n\tilde{B}^m_t\left(u_{m,l}^n(\kappa_1(t))\right)\tilde{\Pi}_l\right\rVert_{2}^2\mathrm{d} t. 
	\end{align*} 
	Here, we used the fact that $\left\lVert \Pi_n\tilde{B}^m_t\left(u_{m,l}^n(\kappa_1(t))\right)\tilde{\Pi}_l\right\rVert_{2}^2$ is constant for $s\in (t,\kappa_2(t))$. Thus by (iv) of step 1, $R_2\to 0$ when $m\to \infty$, $n\geq N$ and $(n,m,l)\in I_\gamma$ for $\gamma \in(0,1)$. Finally, using It\^o-L{\'e}vy's isometry,  the computation for $R_3$ is as follows: 
	\begin{align*} 
		\MoveEqLeft[0] R_3^2 = \left\lvert\mathbb{E}\int_0^T\varphi(t)\int_{t}^{\kappa_2(t)}\int_E \left(\Pi_n\tilde{F}^{m,l}_s\left(u_{m,l}^n(\kappa_1(s)),\xi\right),z\right)_H\tilde{N}(\mathrm{d} s,\mathrm{d}\xi)\mathrm{d} t\right\rvert^2\\ 
		&\leq \left\lVert \varphi\right\rVert_{\mathscr{L}_\mathbb{R}^2}^2\mathbb{E}\int_0^T\left\lvert \int_{t}^{\kappa_2(t)}\int_E \left(\Pi_n\tilde{F}^{m,l}_s\left(u_{m,l}^n(\kappa_1(s)),\xi\right),z\right)_H\tilde{N}(\mathrm{d} s,\mathrm{d}\xi)\right\rvert^2\mathrm{d} t\\ 
		&\leq \left\lVert z\right\rVert^2_H \left\lVert \varphi\right\rVert_{\mathscr{L}_\mathbb{R}^2}^2\mathbb{E}\int_0^T \int_{t}^{\kappa_2(t)}\int_E \left\lVert\Pi_n\tilde{F}^{m,l}_s\left(u_{m,l}^n(\kappa_1(s)),\xi\right)\right\rVert^2_H\nu(\mathrm{d} \xi)\, \mathrm{d} s\,\mathrm{d} t\\&\leq 
		\delta_m\left\lVert z\right\rVert^2_H \left\lVert \varphi\right\rVert_{\mathscr{L}_\mathbb{R}^2(\lambda)}^2\mathbb{E}\int_0^T\int_E \left\lVert \Pi_n\tilde{F}^{m,l}_t\left(u_{m,l}^n(\kappa_1(t)),\xi\right)\right\rVert_{H}^2\nu(\mathrm{d}\xi)\mathrm{d} t
	\end{align*} 
where we used the fact that $\left\lVert\Pi_n\tilde{F}^{m,l}_s\left(u_{m,l}^n(\kappa_1(s)),\xi\right)\right\rVert^2_H$ is constant for $s\in (t,\kappa_2(t))$. This together with (v) of step 1 implies that $R_3\to 0$, when $m\to \infty$, $n\geq N$ and $(n,m,l)\in I_\gamma$ for $\gamma \in(0,1)$. 
	Now we have proven that the limit of the right hand side of equation \eqref{eq:Ji Ri's} is the right hand side of equation \eqref{effect of v infty}. By the fact that $u_{m,l}^n\rightharpoonup \bar{u}_\infty$ we deduce the similar result for the left hand side, so equation \eqref{v infty} is obtained. 
	It remains to prove \eqref{u infty}. Both sides of this equation belong to the Hilbert space $L^2(\Omega;\mathbb{R})$, so it is sufficient to prove that the inner product of both sides with $\psi \in L^2(\Omega;\mathbb{R})$ and $z\in V_N$ are the same. Thus we wish to verify the following equality 
	\begin{align*} 
		\mathbb{E}\left[\psi\left(u_{T\infty},z\right)_H\right]&=\mathbb{E}\left[\psi(\zeta ,z)_H\right]+\mathbb{E}\left[\psi\int_0^T\left\langle 
		a_\infty(s),z\right\rangle \mathrm{d} s\right]\\&\quad+\mathbb{E}\left[\psi\int_0^T\left(b_\infty(s)\mathrm{d} W_s,z\right)_H\right]\\&\quad +\mathbb{E}\left[\psi\int_0^T\int_E\left(f_\infty(s,\xi),z\right)_H\tilde{N}(\mathrm{d} s,\mathrm{d}\xi)\right], 
	\end{align*} 
	for $\psi \in L^2(\Omega;\mathbb{R})$ and $z\in V_N$. 
	Fix $N\in\mathbb{N}$ and $z\in V_N$. By equation \eqref{ex int form}, we get for $n\geq N$ that 
	\begin{align*} 
		(u_{m,l}^n(T),z)_H&=(\zeta ,z)_H\mathbf{1}_{\left\lbrace t>t_0\right\rbrace}+\int_0^{T-\delta_m}\left\langle A_s\left(u_{m,l}^n(s)\right),z\right\rangle\mathrm{d} s\\&\quad +\int_0^{T}\left(z,\Pi_n\tilde{B}^m_s\left(u_{m,l}^n\left(\kappa_1(s)\right)\right)\tilde{\Pi}_l\mathrm{d} W_s\right)_H\\ 
		&\quad+\int_0^{T}\int_E\left(\Pi_n\tilde{F}^{m,l}_s\left(u_{m,l}^n\left(\kappa_1(s)\right)\right),z\right)_H\tilde{N}(\mathrm{d} s,\mathrm{d}\xi) \, . 
	\end{align*} 
	Taking the inner products of both sides and $\psi$, we get for $n\geq N$ 
	\begin{equation}\label{eq:Ji Ri's tilde} 
		\mathbb{E}\left[\psi(u_{m,l}^n(T),z)_H\right]=\mathbb{E}\left[\psi(\zeta ,z)_H\right]+\tilde{J}_1+\tilde{J}_2+\tilde{J}_3-\tilde{R}_1 \, , 
	\end{equation} 
	where 
	\begin{align*} 
		\tilde{J}_1&=\mathbb{E}\left[\psi\int_0^T \left\langle A_s\left(u_{m,l}^n(s)\right),z\right\rangle \mathrm{d} s\right],\\ 
		\tilde{J}_2&=\mathbb{E}\left[\psi\int_0^T\left( \Pi_n\tilde{B}^m_s\left(u_{m,l}^n(\kappa_1(s))\right)\tilde{\Pi}_l\mathrm{d} W_s,z\right)_H\right],\\ 
		\tilde{J}_3&=\mathbb{E}\left[\psi\int_0^T\int_E \left(\Pi_n\tilde{F}^{m,l}_s\left(u_{m,l}^n(\kappa_1(s)),\xi\right),z\right)_H\tilde{N}(\mathrm{d} s,\mathrm{d}\xi)\right],\\ 
		\tilde{R}_1&=\mathbb{E}\left[\psi\int_{T-\delta_m}^T \left\langle A_s\left(u_{m,l}^n(s)\right),z\right\rangle \mathrm{d} s\right].\\ 
	\end{align*} 
	Assume that $\tilde{S}_1$ is the linear operator from $\mathscr{L}_{V^*}^q(\lambda^{1-q},\mathcal{BF})$ to $L^q(\Omega;\mathbb{R})$, defined by 
	\[\tilde{S}_1(g):=\int_0^T \left\langle g(s),z\right\rangle \mathrm{d} s, \] 
	for all $g\in\mathscr{L}_{V^*}^q(\lambda^{1-q},\mathcal{BF})$. $\tilde{S}_1$ is bounded, so it is continuous with respect to the weak topologies. This continuity gives that 
	\[\tilde{J}_1=\mathbb{E}\left[\psi\int_0^T \left\langle A_s\left(u_{m,l}^n(s)\right),z\right\rangle \mathrm{d} s\right]\to\mathbb{E}\left[\psi\int_0^T \left\langle a_\infty(s),z\right\rangle \mathrm{d} s\right].\] 
	For $\tilde{J}_2$, take the linear operator $\tilde{S}_2:\mathscr{L}_{H}^2(\mathcal{BF})\to L^2(\Omega;\mathbb{R})$ as 
	\[\tilde{S}_2(g):=\int_0^T\left( g(s)\mathrm{d} W_s,z\right)_H\] 
	and deduce 
	\begin{align*}
		\tilde{J}_2&=\mathbb{E}\left[\psi\int_0^T\left( \Pi_n\tilde{B}^m_s\left(u_{m,l}^n(\kappa_1(s))\right)\tilde{\Pi}_l\mathrm{d} W_s,z\right)_H\right]\\&\to\mathbb{E}\left[\psi\int_0^T\left( b_\infty(s)\mathrm{d} W_s,z\right)_H\right].\end{align*}
	Similarly, the linear operator
	\[\tilde{S}_3:L^2([0,T]\times\Omega\times E, \mathcal{P}\otimes\mathcal{E},\mathrm{d} t\otimes\mathbb{P}\otimes\nu;H)\to L^2(\Omega;\mathbb{R}),\]\[\tilde{S}_3(g):=\int_0^T\int_E g(s,\xi)\tilde{N}(\mathrm{d} s,\mathrm{d}\xi)\]
	is bounded and continuous with respect to the weak topologies. Therefore
	\begin{align*}
		\tilde{J}_3&=\mathbb{E}\left[\psi\int_0^T\int_E \left(\Pi_n\tilde{F}^{m,l}_s\left(u_{m,l}^n(\kappa_1(s)),\xi\right),z\right)_H\tilde{N}(\mathrm{d} s,\mathrm{d}\xi)\right]\\&\to \mathbb{E}\left[\psi\int_0^T\int_E \left(f_\infty(s,\xi),z\right)_H\tilde{N}(\mathrm{d} s,\mathrm{d}\xi)\right]. 
	\end{align*} 
	It is easy to check that $\tilde{R}_1\to 0$ as $m\to \infty$. By using $u_{m,l}^n(T)\rightharpoonup u_{T\infty}$, equation \eqref{u infty} is obtained. 
	
	\paragraph*{\textbf{Step 3. }}
	\textit{
		Let $(n,m,l)$ be a subsequence which satisfies items \ref{converg of explicit:i}-\ref{converg of explicit:vi} of Step 2. Here $n,m,l$ tend to infinity and $\mathcal{C}(n)\max_{1\leq i\leq m}\int_{t_{i-1}}^{t_i}\lambda(s)\mathrm{d} s$ converges to zero. Then for all $y\in \mathscr{G}_{\mathcal{BF}}$ 
		\[I_y(\zeta ,a_\infty,b_\infty,f_\infty)+\liminf \mathbb{E}\left\lVert u^n_{m,l}(T)\right\rVert_H^2-\mathbb{E}\left\lVert u_{T\infty}\right\rVert_H^2\leq 0.\]} 
	\paragraph*{\text{Proof of Step 3. }}
	Define a bounded linear operator $\mathcal{T}_m$ from $
	\mathcal{H}=L^2([0,T]\times\Omega\times E,\mathcal{BF}\otimes \mathcal{E},\mathrm{d} t\otimes \mathbb{P}\otimes \nu;H)$ or $\mathscr{L}^2_{L_2(U,H)}(\mathcal{BF})$ into themselves as 
	\begin{equation}\label{Tm-define}
		\mathcal{T}_m(Z)(t):=\begin{dcases} 
			\delta_m^{-1}\int_{t_{i-2}}^{t_{i-1}}Z(s)\mathrm{d} s& t_{i-1}<t\leq t_{i},2\leq i\leq m, \\0& otherwise. 
		\end{dcases}
	\end{equation}
	It is easy to check that the induced operator norm of $\mathcal{T}_m$ is equal to $1$. We know from functional analysis that for every $Z\in \mathcal{H}$ or $Z\in \mathscr{L}^2_{L_2(U,H)}(\mathcal{BF})$, when $m\to \infty$, $\mathcal{T}_m(Z)$ converges to $Z$ in $\mathcal{H}$ or $ \mathscr{L}^2_{L_2(U,H)}(\mathcal{BF})$, respectively. Define the operator $\mathcal{S}_l$ on $\mathcal{H}$ as 
	\begin{equation}\label{Sl-define}
		\mathcal{S}_l(Z)(t,\xi):=\begin{dcases} 
			(\nu(E^l_j))^{-1}\int_{E^l_j}Z(t,\eta)\nu(\mathrm{d} \eta)& \xi \in E^l_j, 1\leq j\leq r_l,\\0& otherwise, 
		\end{dcases}
	\end{equation}
	with the convention $0/0:=0$. The operator norm of $\mathcal{S}_l$ is also $1$. Since $\nu\otimes \mathrm{d} t$ is $\sigma$-finite Borel measure, there is a sequence of random functions $Z_n(\omega)\in C([0,T]\times E), \omega\in\Omega$, having a support contained in $[0,T]\times E_n$, such that $Z_n\to Z$ in $\mathcal{H}$ as $n\to\infty$. Let $\varepsilon>0$ be arbitrary. For large enough $n\in \mathbb{N}$, we have
	\[\begin{split}\MoveEqLeft[1]\mathbb{E}\int_0^T\int_E \left\lVert Z(t,\xi)-\mathcal{S}_l(Z)(t,\xi)\right\rVert_H^2\nu(\mathrm{d} \xi)\mathrm{d} t\\&\leq \mathbb{E}\int_0^T\int_{E^n} \left\lVert Z_n(t,\xi)-\mathcal{S}_l(Z_n)(t,\xi)\right\rVert_H^2\nu(\mathrm{d} \xi)\mathrm{d} t+\varepsilon.\end{split}\]
	Since the diameter of $E_j^l, 1\leq j \leq r_l$ is less than $\varepsilon_l$ and $\varepsilon_l\to 0$ as $l\to\infty$, we get by continuity of $Z_n$ that
	\[\begin{aligned}\MoveEqLeft[2]\lim_{l\to\infty}\mathbb{E}\int_0^T\int_E \left\lVert Z(t,\xi)-\mathcal{S}_l(Z)(t,\xi)\right\rVert_H^2\nu(\mathrm{d} \xi)\mathrm{d} t\\&\leq \lim_{l\to\infty}\mathbb{E}\int_0^T\int_{E^n} \left\lVert Z_n(t,\xi)-\mathcal{S}_l(Z_n)(t,\xi)\right\rVert_H^2\nu(\mathrm{d} \xi)\mathrm{d} t+\varepsilon=\varepsilon.\end{aligned}\]
	Since $\varepsilon>0$ is arbitrary, this implies that $\mathcal{S}_l(Z)$ converges to $Z$ in $\mathcal{H}$ as $l\to \infty$.

		We observe that 
		\begin{equation*}
			\begin{split}
				\mathcal{T}_m(\Pi_nB_\cdot\left(u_{m,l}^n(\cdot)\right)\tilde{\Pi}_l)&=\Pi_n\tilde{B}^m_\cdot\left(u_{m,l}^n(\kappa_1(\cdot))\right)\tilde{\Pi}_l
				\\\mathcal{S}_l\circ \mathcal{T}_m(\Pi_nF_\cdot\left(u_{m,l}^n(\cdot),\star\right))&=\Pi_n\tilde{F}^{m,l}_\cdot\left(u_{m,l}^n(\kappa_1(\cdot)), \star\right)
			\end{split}
	\end{equation*}
and for $y\in \mathscr{G}_{\mathcal{BF}}$,
\[\mathcal{T}_m(\Pi_nB_\cdot(y_\cdot))\to B_\cdot (y_\cdot) \qquad \text{in } \mathscr{L}^2_{L_2(U,H)}(\mathcal{BF}),\]
\[\mathcal{S}_l\circ \mathcal{T}_m(\Pi_nF_\cdot(y_\cdot,\star))\to F_\cdot(y_\cdot,\star)\qquad \text{in }\mathcal{H}. \]
	Define 
	\begin{align*} 
		&I_{m,l}^n(y):=\mathbb{E}\int_0^T \Bigg[2\left\langle A_s\left(u_{m,l}^n(s)\right)-A_s(y_s), u_{m,l}^n(s)-y_s\right\rangle\\&+\left\lVert \Pi_n\tilde{B}^m_s\left(u_{m,l}^n(\kappa_1(s))\right)\tilde{\Pi}_l-\mathcal{T}_m(\Pi_nB_\cdot(y_\cdot))(s)\right\rVert_2^2\\&+\int_E\left\lVert \Pi_n\tilde{F}^{m,l}_s\left(u_{m,l}^n(\kappa_1(s)),\xi\right)-\mathcal{S}_l\circ \mathcal{T}_m(\Pi_nF_\cdot(y_\cdot,\star))(s,\xi)\right\rVert_H^2\nu(\mathrm{d} \xi)\Bigg]\mathrm{d} s. 
	\end{align*} 
Since the operator norms of $\mathcal{T}_m$ and $\mathcal{S}_l$ equal one, we have
	\begin{align*} 
		I_{m,l}^n(y)&\leq\mathbb{E}\int_0^T \Bigg[2\left\langle A_s\left(u_{m,l}^n(s)\right)-A_s(y_s), u_{m,l}^n(s)-y_s\right\rangle\\&\quad+\left\lVert \Pi_n B_s\left(u_{m,l}^n(s)\right)\tilde{\Pi}_l-\Pi_nB_s(y_s)\right\rVert_2^2\\&\quad+\int_E\left\lVert \Pi_nF_s\left(u_{m,l}^n(s),\xi\right)-\Pi_nF_s(y_s,\xi)\right\rVert_H^2\nu(\mathrm{d} \xi)\Bigg]\mathrm{d} s \, . 
	\end{align*}
	Using the monotonicity condition \ref{C1}, it is obvious to see that $I_{m,l}^n(y)\leq 0$. By taking the sum of inequality \eqref{before15} over $2\leq i\leq m$, we get
	\begin{align*} 
		\MoveEqLeft[0]\mathbb{E}\left\lVert u^n_{m,l}(T)\right\rVert^2_H\\&\leq \mathbb{E}\left\lVert \zeta \right\rVert_H^2+\left(\max_{1\leq i\leq m}\int_{t_{i-1}}^{t_i}\lambda(s)\mathrm{d} s\right) \mathcal{C}(n)\mathbb{E}\int_0^{T-\delta_m}\left\lVert A_s\left(u_{m,l}^n(s)\right)\right\rVert_{V^*}^2\lambda(s)^{-1}\mathrm{d} s\\&\quad+\mathbb{E}\int_0^T\Bigg[ 2\left\langle u_{m,l}^n(s),A_s\left(u_{m,l}^n(s)\right)\right\rangle+\left\lVert \Pi_n \tilde{B}^m_s\left(u_{m,l}^n(\kappa_1(s))\right)\tilde{\Pi}_l\right\rVert_2^2 \\& \quad+\int_E\left\lVert \Pi_n \tilde{F}^{m,l}_s\left(u_{m,l}^n(\kappa_1(s)),\xi\right)\right\rVert_H^2\nu(\mathrm{d} \xi)\Bigg]\mathrm{d} s\\&\quad - \mathbb{E}\int_{T-\delta_m}^T2\left\langle u_{m,l}^n(s),A_s\left(u_{m,l}^n(s)\right)\right\rangle \mathrm{d} s\, . 
	\end{align*} 
	Thus 
	\begin{equation}\label{I^n_m,l} 
		\begin{split}
			0\geq I_{m,l}^n(y)&\geq \mathbb{E}\left\lVert u^n_{m,l}(T)\right\rVert_H^2-\mathbb{E}\left\lVert \zeta \right\rVert_H^2+\mathbb{E}\int_0^T 2\left\langle A_s(y_s),y_s\right\rangle\mathrm{d} s\\&\quad+J_1+J_2-2(J_3+J_4+J_5+J_6)-R_1+2R_2 \, , 
		\end{split} 
	\end{equation}
	where 
	\begin{align*} 
		J_1&:=\mathbb{E}\int_0^T\left\lVert \mathcal{T}_m(\Pi_n B_\cdot(y_\cdot))(s)\right\rVert_2^2\mathrm{d} s\to \mathbb{E}\int_0^T\left\lVert B_s(y_s)\right\rVert_2^2\mathrm{d} s, 
		\\J_2&:=\mathbb{E}\int_0^T\int_E\left\lVert \mathcal{S}_l\circ \mathcal{T}_m\left(\Pi_n F_\cdot(y_\cdot,\star\right)(s,\xi)\right\rVert_H^2\nu(\mathrm{d} \xi)\mathrm{d} s\\&\to \mathbb{E}\int_0^T\int_E\left\lVert F_s(y_s,\xi)\right\rVert_H^2\nu(\mathrm{d} \xi)\mathrm{d} s,\\ 
		J_3&:=\mathbb{E}\int_0^T\left\langle A_s\left(u_{m,l}^n(s)\right),y_s\right\rangle\mathrm{d} s\to \mathbb{E}\int_0^T \left\langle a_\infty(s),y_s\right\rangle\mathrm{d} s,\\ 
		J_4&:=\mathbb{E}\int_0^T\left\langle A_s(y_s),u^n_{m,l}(s)\right\rangle\mathrm{d} s\to \mathbb{E}\int_0^T\left\langle A_s(y_s),\bar{u}_\infty(s)\right\rangle\mathrm{d} s,\\ 
		J_5&:=\mathbb{E}\int_0^T\left\langle \mathcal{T}_m(\Pi_n B_\cdot(y_\cdot))(s),\Pi_n\tilde{B}^m_s\left(u^n_{m,l}(\kappa_1(s))\right)\tilde{\Pi}_l\right\rangle_2\mathrm{d} s\\&\to \mathbb{E}\int_0^T\left\langle B_s(y_s),b_\infty(s)\right\rangle_2\mathrm{d} s
		,\\
		J_6&:=\mathbb{E}\int_0^T\int_E\left(\mathcal{S}_l\circ \mathcal{T}_m\left(\Pi_n F_\cdot(y_\cdot,\star)\right)(s,\xi),\Pi_n \tilde{F}^{m,l}_s\left(u^n_{m,l}(\kappa_1(s)),\xi\right)\right)_H\nu(\mathrm{d} \xi)\mathrm{d} s\\&\to \mathbb{E}\int_0^T\int_E\left( F_s(y_s,\xi),f_\infty\left(s,\xi\right)\right)_H\nu(\mathrm{d} \xi)\mathrm{d} s,
	\end{align*}
	and
	\begin{align*} 
		R_1&:=\left(\max_{1\leq i\leq m}\int_{t_{i-1}}^{t_i}\lambda(s)\mathrm{d} s\right) \mathcal{C}(n)\mathbb{E}\int_0^{T-\delta_m}\left\lVert A_s\left(u_{m,l}^n(s)\right)\right\rVert_{V^*}^2\lambda(s)^{-1}\mathrm{d} s\to 0,\\ 
		R_2&:=\mathbb{E}\int_{T-\delta_m}^T\left\langle u_{m,l}^n(s),A_s\left(u_{m,l}^n(s)\right)\right\rangle \mathrm{d} s. 
	\end{align*} 
	Now it remains to identify the limit of $R_2$. By parts \ref{boundedness:i} and \ref{boundedness:iii} of Step 1  and Cauchy-Schwartz's inequality, we get that 
	\begin{align*} 
		R_2^2&\leq\left(\mathbb{E}\int_{T-\delta_m}^T\left\lVert u^n_{m,l}(s)\right\rVert_H^2\lambda(s)\mathrm{d} s\right)\left(\mathbb{E}\int_0^T\left\lVert \Pi_n A_s\left(u^n_{m,l}(s)\right)\right\rVert_H^2\lambda^{-1}(s)\mathrm{d} s\right)\\ &\leq \left(\max_{1\leq i\leq m}\int_{t_{i-1}}^{t_i}\lambda(s)\mathrm{d} s\right)\mathcal{C}(n)\sup_{0\leq i\leq m}\mathbb{E}\left\lVert u^n_{m,l}(t_{i})\right\rVert_H^2\left\lVert A_\cdot\left(u^n_{m,l}(\cdot)\right)\right\rVert_{\mathscr{L}^2_{V^*}(\lambda^{-1})}^2\\&\to 0 \, . 
	\end{align*} 
	%
	Now the limits of all terms in the inequality \eqref{I^n_m,l}, except $\mathbb{E}\left\lVert u^n_{m,l}(T)\right\rVert_H^2$, have been identified. Since $u^n_{m,l}(T)\rightharpoonup u_{T\infty}$, we can write 
	\begin{equation}\label{d-explicit}
		d:=\liminf \mathbb{E}\left\lVert u^n_{m,l}(T)\right\rVert_H^2-\mathbb{E}\left\lVert u_{T\infty}\right\rVert_H^2\geq 0.
	\end{equation}
	By Theorem \ref{theorem1of1}, $\bar{u}_\infty$ has an adapted c\`adl\`ag $H$-valued modification, denoted by $u_\infty$, such that almost surely for all $t\in[0,T]$ 
	\[
	u_\infty(t)=\zeta +\int_0^t a_\infty(s)\mathrm{d} s+\int_0^t b_\infty (s)\mathrm{d} W_s+\int_0^t\int_E f_\infty(s,\xi)\tilde{N}(\mathrm{d} s, \mathrm{d} \xi). 
	\] 
	This equation for $t=T$, together with the equation \eqref{u infty}, implies that 
	$u_\infty (T)=u_{T\infty}$ a.s. 
	By It\^o's formula (see Theorem \ref{theorem1of1}), we have that 
	\begin{equation}\label{uTinfty-explicit}
	\begin{split}
		&\mathbb{E}\left\lVert u_{T\infty}\right\rVert_H^2 \\&= \mathbb{E}\left\lVert \zeta \right\rVert_H^2 +\mathbb{E}\int_0^T\left[ 2\left\langle a_\infty(s),\bar{u}_\infty(s)\right\rangle+\left\lVert b_\infty(s)\right\rVert_2^2+\int_E\left\lVert f_\infty(s,\xi)\right\rVert_H^2\nu(\mathrm{d} \xi)\right]\mathrm{d} s. 
		\end{split}
	\end{equation} 
	Therefore, by taking limit inferior of inequality \eqref{I^n_m,l} and using \eqref{d-explicit} and \eqref{uTinfty-explicit}, we get
		\begin{align*} 
			0&\geq d +\mathbb{E}\int_0^T\left[ 2\left\langle a_\infty(s),\bar{u}_\infty(s)\right\rangle+\left\lVert b_\infty(s)\right\rVert_2^2+\int_E\left\lVert f_\infty(s,\xi)\right\rVert_H^2\nu(\mathrm{d} \xi)\right]\mathrm{d} s\\&+\mathbb{E}\int_0^T \left[2\left\langle A_s(y_s),y_s\right\rangle+\left\lVert B_s(y_s)\right\rVert_2^2+\int_E\left\lVert F_s(y_s,\xi)\right\rVert_H^2\nu(\mathrm{d} \xi)\right]\mathrm{d} s\\& -2\, \mathbb{E}\int_0^T \left[\left\langle a_\infty(s),y_s\right\rangle+\left\langle A_s(y_s),\bar{u}_\infty(s)\right\rangle\right]\mathrm{d} s \\&-2\, \mathbb{E}\int_0^T \left[\left\langle B_s(y_s),b_\infty(s)\right\rangle_2+\int_E\left( F_s(y_s,\xi),f_\infty\left(s,\xi\right)\right)_H\nu(\mathrm{d} \xi)\right]\mathrm{d} s \, , \\&=
			d+I_y(\zeta ,a_\infty,b_\infty,f_\infty). 
	\end{align*}
	Hence by Theorem \ref{characterize}, $u_\infty$ must be a solution of equation \eqref{equation1} 
	which is unique. Setting $y=u_\infty$, we get $d=0$ which implies that a subsequence of 
	$u^n_{m,l}(T)$ converges strongly in $L^2\left(\Omega,\mathbb{P};H\right)$ to $u_\infty(T)$. Since 
	there exist unique limits for any convergent subsequences of $u^n_{m,l}$ and $u^n_{m,l}(T)$, and on the other hand 
		each subsequence of these sequences has a convergent subsequence, the whole sequences $u^n_{m,l}$ and $u^n_{m,l}(T)$ converge. Thus the proof of 
	theorem is complete. 
\end{proof} 

\subsection{Convergence of the Implicit Schemes} 

The poof of Theorem \ref{implicitexistance} is based on \cite[Proposition 3.4]{gyongy2005discretization}. 
\begin{proposition}\label{prop3} 
	\cite{gyongy2005discretization} Let $D:V\to V^*$ be such that: 
	\begin{enumerate}[label=(\roman{*})] 
		\item $D$ is monotone, i.e. for every $x,y\in V$, $\left\langle D(x)-D(y),x-y\right\rangle\geq 0$. 
		\item $D$ is hemicontinuous, i.e. $\lim_{\varepsilon\to 0}\left\langle D(x+\varepsilon y),z\right\rangle=\left\langle D(x),z\right\rangle$ for every $x,y,z\in V$. 
		\item $D$ satisfies the growth condition, i.e. there exists $K>0$ such that for every $x\in V$, 
		\[\left\lVert D(x)\right\rVert_{V^*}\leq K(1+\left\lVert x\right\rVert_V^{p-1}).\] 
		\item $D$ is coercive, i.e. there exist constants $C_1>0$ and $C_2\geq 0 $ such that 
		\[\left\langle D(x),x\right\rangle\geq C_1\left\lVert x\right\rVert_V^p-C_2,\quad \forall x\in V\] 
	\end{enumerate} 
	Then for every $y\in V^*$, there exists $x\in V$ such that $D(x)=y$ and 
	\[\left\lVert x\right\rVert_V^p\leq \frac{C_1+2C_2}{C_1}+\frac{1}{C_1^2}\left\lVert y\right\rVert_{V^*}^2.\] 
	If there exists a positive constant $C_3$ such that 
	\begin{equation}\label{dissipative D}
		\left\langle D(x_1)-D(x_2),x_1-x_2\right\rangle\geq C_3\left\lVert x_1-x_2\right\rVert_{V^*}^2,\quad x_1,x_2\in V,
	\end{equation}
	then for any $y\in V^*$, the equation $D(x)=y$ has a unique solution $x\in V$. 
\end{proposition} 
\begin{proof}[Proof of Theorem \ref{implicitexistance}] 
	Let $Id:V\to V$ and $Id_n:V_n\to V_n$ be the identity maps. It is easy to check that for large enough $m$, the operators
	\[
	D=Id-\int_{t_{i-1}}^{t_{i}} A_s(x)\mathrm{d} s\, , 
	\] 
	\[
	D_n=Id_n-\int_{t_{i-1}}^{t_{i}} \Pi_n A_s(x)\mathrm{d} s\, , 
	\] 
	have the properties (i)-(iv) stated in Proposition \ref{prop3} and satisfy \eqref{dissipative D}. Using Proposition \ref{prop3}, each of equations $D(x)=y$ and $D_n(x)=y$ has unique solution. According to \eqref{dissipative D}, $D^{-1}$ and $D_n^{-1}$ are continuous and hence they are measurable. 
	So the existence of unique $\mathcal{F}_{t_i}$-measurable solutions to \eqref{implicitnotspace} and \eqref{implicit} follows inductively. 
	We use induction on $i$ to prove that $\mathbb{E}\left\lVert u^{n,m,l}(t_i)\right\rVert_V^p< \infty$. The proof of $\mathbb{E}\left\lVert u^{m,l}(t_i)\right\rVert_V^p<\infty$ is the same. According to the previous proposition, we get that 
	\[
	\left\lVert u^{n,m,l}(t_{i})\right\rVert_V^p\leq C'_1+C'_2\left\lVert  y_{i}\right\rVert_{V^*}^2 \, . 
	\] 
	where 
	\begin{align*} 
	y_{i}&:=u^{n,m,l}(t_{i-1})+\Pi_n\tilde{B}^m_{t_{i}}\left(u^{n,m,l}(t_{i-1})\right)\left(W^l_{t_i}-W^l_{t_{i-1}}\right)\\&\quad+\int_E\Pi_n \tilde{F}^{m,l}_{t_i}\left(u^{n,m,l}(t_{i-1}),\xi\right)\tilde{N}\left((t_{i-1},t_i],\mathrm{d} \xi\right). 
	\end{align*} 
For $i=1$, $y_1=\Pi_n\zeta$ and $\mathbb{E}\left\lVert y_1\right\rVert_{V^*}^2\leq C\mathbb{E}\left\lVert \zeta\right\rVert_{H}^2<\infty$, so we get $\mathbb{E}\left\lVert u^{n,m,l}(t_1) \right\rVert^p_V<\infty$. Let us estimate $\mathbb{E}\left\lVert y_i\right\rVert_{V^*}^2$ for $i\geq 2$. Using Proposition 
	\ref{proposition1} and $p\geq 2$, we have 
	\begin{align*} 
		\MoveEqLeft[1.3]\mathbb{E}\left\lVert y_i\right\rVert_{V^*}^2\leq C\mathbb{E}\left\lVert y_i\right\rVert_{H}^2\\ &\leq C\mathbb{E}\left\lVert u^{n,m,l}(t_{i-1})\right\rVert_H^2+C\delta_m \mathbb{E}\left\lVert \Pi_n\tilde{B}^m_{t_i}\left(u^{n,m,l}(t_{i-1})\right)\right\rVert_2^2+\\&\quad+C\delta_m\int_E\mathbb{E}\left\lVert \Pi_n\tilde{F}_{t_i}^{m,l}\left(u^{n,m,l}(t_{i-1}),\xi\right)\right\rVert_H^2\nu(\mathrm{d} \xi)\\ &\leq 
		C\mathbb{E}\left\lVert u^{n,m,l}(t_{i-1})\right\rVert_H^2+C\delta_m^{-1}\mathbb{E}\left\lVert \int_{t_{i-2}}^{t_{i-1}}\Pi_nB_s\left(u^{n,m,l}(t_{i-1})\right)\mathrm{d} s\right\rVert_2^2 \\&+C\delta_m^{-1}\sum_{1\leq j\leq r_l}\left(\nu(E^l_j)\right)^{-1}\mathbb{E}\left\lVert \int_{t_{i-2}}^{t_{i-1}}\int_{E^l_j} \Pi_n F_s\left(u^{n,m,l}(t_i),\xi\right)\nu(\mathrm{d} \xi)\mathrm{d} s\right\rVert_H^2\\& \leq C\mathbb{E}\left\lVert u^{n,m,l}(t_{i-1})\right\rVert_H^2\\ &+C\mathbb{E}\int_{t_{i-2}}^{t_{i-1}}\left[\left\lVert B_s\left( u^{n,m,l}(t_{i-1})\right)\right\rVert_2^2\mathrm{d} s+\int_E \left\lVert F_s\left(u^{n,m,l}(t_{i-1}),\xi\right)\right\rVert_H^2\nu(\mathrm{d} \xi)\right]\mathrm{d} s\\ &\leq C\mathbb{E}\left\lVert u^{n,m,l}(t_{i-1})\right\rVert_H^2 \\&+C\mathbb{E}\int_{t_{i-2}}^{t_{i-1}}\left[ 6\left(\alpha^{1/q}+1/p\right) \left\lVert u^{n,m,l}(t_{i-1})\right\rVert_V^p\lambda(s)+\frac{16}{q}K_2(s)+2K_1(s)\right]\mathrm{d} s\\&\leq C\left(1+\mathbb{E}\left\lVert u^{n,m,l}(t_{i-1})\right\rVert_V^p\right). 
	\end{align*} 
	So, if $\mathbb{E}\left\lVert u^{n,m,l}(t_{i-1})\right\rVert_V^p<\infty$, then $\mathbb{E}\left\lVert u^{n,m,l}(t_{i})\right\rVert_V^p<\infty$ too, and the proof of Theorem \ref{implicitexistance} is complete. 
\end{proof} 
\begin{proof}[Proof of Theorem \ref{implicitconvergence}]
	The integral form of the implicit scheme is 
	\begin{equation}\label{int form of imp} 
		\begin{split}
			u^{n,m,l}(t)&=\Pi_n \zeta +\int_0^{\kappa_2(t)}\Pi_n A_s\left(u^{n,m,l}(s)\right)\mathrm{d} s\\&\quad+\int_0^{\kappa_2(t)}\Pi_n\tilde{B}^m_s\left(u^{n,m,l}\left(\kappa_1(s)\right)\right)\mathrm{d} W^l_s\\&\quad+\int_0^{\kappa_2(t)}\int_E\Pi_n\tilde{F}^{m,l}_s\left(u^{n,m,l}\left(\kappa_1(s)\right),\xi\right)\tilde{N}(\mathrm{d} s,\mathrm{d} \xi). 
		\end{split} 
	\end{equation}
	\paragraph*{\textbf{Step 1. }} 
	\textit{There exist a natural number $m_0$ and a constant $L>0$ such that for every $m\geq m_0$ and $n,l\geq 1$ the value of 
		\begin{equation}\label{impboundof unml}\begin{split} 
				\MoveEqLeft\left\lVert u^{n,m,l}\right\rVert_{L^\infty([0,T],\mathrm{d} t; L^2(\Omega;H))}+\left\lVert u^{n,m,l}\right\rVert_{\mathscr{L}^p_V(\lambda)}+\left\lVert A_\cdot\left(u^{n,m,l}(\cdot)\right)\right\rVert_{\mathscr{L}^q_{V^*}(\lambda^{1-q})}\\&\quad+\left\lVert \Pi_n \tilde{B}^m_\cdot\left(u^{n,m,l}\left(\kappa_1(\cdot)\right)\right)\tilde{\Pi}_l\right\rVert_{\mathscr{L}^2_{L_2(U,H)}}\\&\quad+\left\lVert \Pi_n \tilde{F}^{m,l}_\cdot\left(u^{n,m,l}\left(\kappa_1(\cdot)\right),*\right)\right\rVert_{L^2([0,T]\times\Omega\times E;H)} \end{split}
		\end{equation} 
		is bounded by $L$. The same is true for $u^{m,l}$. }
	\paragraph*{\text{Proof of Step 1. }}
	We prove this step only for $u^{n,m,l}$. The proof for $u^{m,l}$ can be done similarly. Since $\mathcal{F}_{t_{i-1}}$, $W_{t_{i}}-W_{t_{i-1}}$, and $\tilde{N}((t_{i-1},t_i], \mathrm{d} \xi)$ are independent, we get from \eqref{implicit} that
	\begin{equation}\label{ghablan(4)imp}
		\begin{split}
			\MoveEqLeft[2]\mathbb{E}\left\lVert u^{n,m,l}(t_i)\right\rVert_H^2- 2\delta_m\mathbb{E}\left\langle \Pi_n A^m_{t_i}\left(u^{n,m,l}(t_i)\right),u^{n,m,l}(t_i)\right\rangle\\&\leq \mathbb{E}\left\lVert u^{n,m,l}(t_{i})-\delta_m\Pi_n A_{t_i}^m\left(u^{n,m,l}
			(t_{i})\right)\right\rVert_H^2\\&\leq \mathbb{E}\left\lVert u^{n,m,l}(t_{i-1})\right\rVert_H^2+\delta_m\mathbb{E} \left\lVert \Pi_n \tilde{B}^m_{t_i}\left(u^{n,m,l}(t_{i-1})\right)\tilde{\Pi}_l\right\rVert_2^2\\&\quad+\delta_m \int_E\mathbb{E}\left\lVert \Pi_n\tilde{F}^{m,l}_{t_i}\left(u^{n,m,l}(t_{i-1}),\xi\right)\right\rVert_H^2\nu(\mathrm{d} \xi)
		\end{split}
	\end{equation}
	So
	\begin{equation}
		\begin{split}
			\mathbb{E}\left\lVert u^{n,m,l}(t_i)\right\rVert_H^2&\leq \mathbb{E}\left\lVert u^{n,m,l}(t_{i-1})\right\rVert_H^2+ 2\mathbb{E}\int_{t_{i-1}}^{t_i}\left\langle u^{n,m,l}(t_i),A_s\left(u^{n,m,l}(t_i)\right)\right\rangle\mathrm{d} s\\&\quad +\mathbb{E}\int_{t_{i-2}}^{t_{i-1}}\Bigg[\left\lVert \Pi_nB_s\left(u^{n,m,l}(t_{i-1})\right)\tilde{\Pi}_l\right\rVert_2^2\\&\qquad \qquad+\int_E\left\lVert\Pi_n F_s\left(u^{n,m,l}(t_{i-1}),\xi\right)\right\rVert_H^2\nu(\mathrm{d} \xi)\Bigg]\mathrm{d} s.
		\end{split}
	\end{equation} 
	Summing up the above inequalities with respect to $i$ and using the coercivity condition we obtain that 
	\begin{align*} 
		\mathbb{E}\left\lVert u^{n,m,l}(t_i)\right\rVert_H^2&\leq \mathbb{E}\left\lVert \zeta \right\rVert_H^2+\mathbb{E}\int_0^{t_i}\bigg[ 2\left\langle u^{n,m,l}(s),A_s\left(u^{n,m,l}(s)\right)\right\rangle\\&\quad+\left\lVert \Pi_n B_s\left(u^{n,m,l}(s)\right)\tilde{\Pi}_l\right\rVert_2^2+\left\lVert \Pi_n F_s\left(u^{n,m,l}(s)\right)\right\rVert_H^2\bigg]\mathrm{d} s\\&\leq \mathbb{E}\left\lVert \zeta \right\rVert_H^2+\mathbb{E}\int_0^{t_i}\bigg[ -\lambda(s)\left\lVert u^{n,m,l}(s)\right\rVert_V^p\\&\quad +\mu(s)\left\lVert u^{n,m,l}(s)\right\rVert_H^2+K_1(s)\bigg]\mathrm{d} s \, . 
	\end{align*} 
	It is evident that there exists a natural number $m_1$ such that for every $m\geq m_1$: 
	\[
	\sup_{1\leq i\leq m}\int_{t_{i-1}}^{t_i}\mu(s)\mathrm{d} s\leq \frac{1}{2}\, . 
	\] 
	So by taking $\alpha_i:=\int_{t_{i-1}}^{t_i}\mu(s)\mathrm{d} s$ for $1\leq i\leq m$, one obtains that 
	\begin{equation} \label{ghablan 16} 
		\frac{1}{2}\mathbb{E}\left\lVert u^{n,m,l}(t_i)\right\rVert_H^2+\mathbb{E}\int_0^{t_i}\lambda(s)\left\lVert u^{n,m,l}(s)\right\rVert_V^p\mathrm{d} s\leq C+\sum_{k=1}^{i-1}\alpha_k\mathbb{E}\left\lVert u^{n,m,l}(t_k)\right\rVert_H^2.
	\end{equation} 
	By using induction on $i$, it is easy to check that 
	\[
	\mathbb{E}\left\lVert u^{n,m,l}(t_i)\right\rVert_H^2\leq 2C(1+2\alpha_1)(1+2\alpha_2)\cdots(1+2\alpha_{i-1}). 
	\] 
	So 
	\[
	\mathbb{E}\left\lVert u^{n,m,l}(t_i)\right\rVert_H^2\leq 2C(1+2\alpha_1)\cdots(1+2\alpha_{m})\leq 2C\left(1+\frac{2\int_0^T\mu(s)\mathrm{d} s}{m}\right)^m \, , 
	\] 
	and this inequality implies that 
	\[
	\sup_{n,l\geq 1,m\geq m_1}\sup_{s\in[0,T]}\mathbb{E}\left\lVert u^{n,m}(s)\right\rVert_H^2 
	< \infty \, . 
	\] 
	Therefore by inequality \eqref{ghablan 16}, the boundedness of $u^{n,m,l}$ in 
	$\mathscr{L}^p_V(\lambda)$ follows. The boundedness of $A_\cdot\left(u^{n,m,l}(\cdot)\right)$ in 
	$\mathscr{L}^q_{V^*}(\lambda^{1-q})$ is obtained from the growth condition \ref{C3}. The 
	boundedness of the fourth and the fifth summand in relation \eqref{impboundof unml} can be 
	obtained in the same way as the proof of parts \ref{boundedness:iv} and \ref{boundedness:v} 
	of Step 1 in the proof of Theorem \ref{explicitconvergence}. 
	
	\paragraph*{\textbf{Step 2. }}
	\textit{ Let $(n,m,l)$ be a sequence such that $m,n$ and $l$ converge to infinity. Then it contains a subsequence, denoted also by $(n,m,l)$, such that 
		\begin{enumerate}[label=(\roman{*})] 
			\item \label{converg of implicit:i}$u^{n,m,l}$ converges weakly in $\mathscr{L}^p_V(\lambda)$ to some progressively measurable process $\bar{u}_\infty$, 
			\item \label{converg of implicit:ii}$u^{n,m,l}(T)$ converges weakly in $L^2(\Omega;H)$ to some random variable $u_{T\infty}$, 
			\item \label{converg of implicit:iii}$A_\cdot \left(u^{n,m,l}(\cdot)\right)$ converges weakly in $\mathscr{L}_{V^*}^q(\lambda^{1-q})$ to some progressively measurable process $a_\infty$, 
			\item \label{converg of implicit:iv}$\Pi_n\tilde{B}_\cdot^m\left(u^{n,m,l}(\kappa_1(\cdot))\right)\tilde{\Pi}_l$ converges weakly in $\mathscr{L}_{L_2(U,H)}^2(\mathcal{BF})$ to some process $b_\infty$, 
			\item \label{converg of implicit:v}$\Pi_n\tilde{F}_\cdot^{m,l}\left(u^{n,m,l}(\kappa_1(\cdot)),*\right)$ converges weakly in $L^2([0,T]\times\Omega\times E, \mathcal{P}\otimes\mathcal{E},\mathrm{d} t\otimes\mathbb{P}\otimes\nu;H)$ to some process $f_\infty$, 
			\item \label{converg of implicit:vi}$(\zeta ,a_\infty,b_\infty,f_\infty)\in \mathcal{A}$ and for all $z\in V$, and $\mathrm{d} t\otimes \mathbb{P}$-almost all $(t,\omega)$ we have 
			\begin{equation}
				\label{v inftyimp}
				\begin{split} 
					(\bar{u}_\infty(t),z)_H&=(\zeta ,z)_H+\int_0^t \langle a_{\infty}(s),z\rangle \mathrm{d} s+\int_0^t \left(b_{\infty}(s)\mathrm{d} W_s,z\right)_H\\&\quad+\int_0^t\int_E(f_{\infty}(s,\xi),z)_H\tilde{N}(\mathrm{d} s,\mathrm{d}\xi) 
				\end{split} 
			\end{equation}
			and for all $z\in V$, almost surely 
			\begin{equation}\label{u inftyimp} 
				\begin{split}
					(u_{T\infty},z)_H&=(\zeta ,z)_H+\int_0^T\langle a_\infty(s),z\rangle \mathrm{d} s+\int_0^T\left(z,b_\infty(s)\mathrm{d} W_s\right)_H\\&\quad+\int_0^T\int_E(f_\infty(s,\xi),z)_H\tilde{N}(\mathrm{d} s,\mathrm{d}\xi)\, . 
				\end{split}
			\end{equation} 
	\end{enumerate} }
	
	\paragraph*{\text{Proof of Step 2. }}
	The convergences in \ref{converg of implicit:i}-\ref{converg of implicit:v} directly follow from 
	\eqref{impboundof unml}. The fact that $\bar{u}_\infty$ and $a_\infty$ are progressively measurable 
	can be shown similarly to the proof of the corresponding statements in Step 2 of the proof of 
	Theorem \ref{explicitconvergence}. It remains to prove \ref{converg of implicit:vi}. Fix 
	$N\in\mathbb{N}$. It is sufficient to verify \ref{converg of implicit:vi} for $z\in V_N$, because 
	$\bigcup_{N=1}^\infty V_N$ is dense in $V$. To verify \eqref{v inftyimp}, it is sufficient to prove 
	that for any $\varphi\in \mathscr{L}_\mathbb{R}^\infty$, 
	\begin{equation} 
		\label{effect of v inftyimp} 
		\begin{split} 
			\MoveEqLeft\mathbb{E}\int_0^T \left(\bar{u}_\infty(t),z\right)_H\varphi(t)\mathrm{d} t\\&=\mathbb{E}\int_0^T\left(\zeta ,z\right)_H\varphi(t)\mathrm{d} t+\mathbb{E}\int_0^T\left(\int_0^t\left\langle a_\infty(s),z\right\rangle \mathrm{d} s\right)\varphi(t)\mathrm{d} t\\&\quad +\mathbb{E}\int_0^T\left(\int_0^t\left(b_\infty(s)\mathrm{d} W_s,z\right)_H\right)\varphi(t)\mathrm{d} t\\&\quad +\mathbb{E}\int_0^T\left(\int_0^t\int_E \left(f_\infty(s,\xi),z\right)_H\tilde{N}(\mathrm{d} s,\mathrm{d}\xi)\right)\varphi(t)\mathrm{d} t. 
		\end{split} 
	\end{equation} 
	Note that the integral form of the implicit scheme \eqref{int form of imp} for $z\in V_N$ and $n\geq N$, yields 
	\begin{align*} 
		(u^{n,m,l}(t),z)_H&=(\zeta ,z)_H+\int_0^{\kappa_2(t)}\left\langle A_s\left(u^{n,m,l}(s)\right),z\right\rangle\mathrm{d} s\\&\quad +\int_0^{\kappa_2(t)}\left(z,\Pi_n\tilde{B}^m_s\left(u^{n,m,l}\left(\kappa_1(s)\right)\right)\mathrm{d} W^l_s\right)_H\\ 
		&\quad+\int_0^{\kappa_2(t)}\int_E\left(\Pi_n\tilde{F}^{m,l}_s\left(u^{n,m,l}\left(\kappa_1(s)\right)\right),z\right)_H\tilde{N}(\mathrm{d} s,\mathrm{d}\xi) \, . 
	\end{align*} 
	Taking the inner products of both sides and $\varphi$, we get for $n\geq N$, 
	\begin{equation}\label{eq:Ji Ri's imp}
		\begin{split} 
			\mathbb{E}\int_0^T\left( u^{n,m,l}(t),z\right)_H\varphi(t)\mathrm{d} t&=\mathbb{E}\int_0^T\left(\zeta ,z\right)_H\varphi(t)\mathrm{d} t+J_1+J_2+J_3\\&\quad-R_1-R_2-R_3 \, , 
		\end{split}
	\end{equation} 
	where 
	\begin{align*} 
		J_1&=\mathbb{E}\int_0^T\varphi(t)\left(\int_0^t \left\langle A_s\left(u^{n,m,l}(s)\right),z\right\rangle \mathrm{d} s\right)\mathrm{d} t\\&\to \mathbb{E}\int_0^T\varphi(t)\left(\int_0^t \left\langle a_\infty(s),z\right\rangle \mathrm{d} s\right)\mathrm{d} t ,\\ 
		J_2&=\mathbb{E}\int_0^T\varphi(t)\int_0^t\left( \Pi_n\tilde{B}^m_s\left(u^{n,m,l}(\kappa_1(s))\right)\tilde{\Pi}_l\mathrm{d} W_s,z\right)_H\mathrm{d} t\\&\to\mathbb{E}\int_0^T\varphi(t)\int_0^t\left( b_\infty(s)\mathrm{d} W_s,z\right)_H\mathrm{d} t,\\ 
		J_3&=\mathbb{E}\int_0^T\varphi(t)\left(\int_0^t\int_E \left(\Pi_n\tilde{F}^{m,l}_s\left(u^{n,m,l}(\kappa_1(s)),\xi\right),z\right)_H\tilde{N}(\mathrm{d} s,\mathrm{d}\xi)\right)\mathrm{d} t\\&\to\mathbb{E}\int_0^T\varphi(t)\left(\int_0^t\int_E \left(f_\infty(s,\xi),z\right)_H\tilde{N}(\mathrm{d} s,\mathrm{d}\xi)\right)\mathrm{d} t,\end{align*}
	and
	\begin{align*}
		R_1&=\mathbb{E}\int_0^T\varphi(t)\left(\int_t^{\kappa_2(t)} \left\langle A_s\left(u^{n,m,l}(s)\right),z\right\rangle \mathrm{d} s\right)\mathrm{d} t,\\ 
		R_2&=\mathbb{E}\int_0^T\varphi(t)\int_{t}^{\kappa_2(t)} \left(\Pi_n\tilde{B}^m_s\left(u^{n,m,l}(\kappa_1(s))\right)\tilde{\Pi}_l\mathrm{d} W_s,z\right)_H\mathrm{d} t,\\ 
		R_3&=\mathbb{E}\int_0^T\varphi(t)\int_{t}^{\kappa_2(t)}\int_E \left(\Pi_n\tilde{F}^{m,l}_s\left(u^{n,m,l}(\kappa_1(s)),\xi\right),z\right)_H\tilde{N}(\mathrm{d} s,\mathrm{d}\xi)\mathrm{d} t.\\ 
	\end{align*} 
	The limits of $J_i$, $i=1,2,3$ can be obtained similar to the limits of $J_i, i=1,2,3$ of \eqref{eq:Ji Ri's}. Now we wish to prove that the "$R_i$"s tend to zero. H\"older's inequality yields 
	\begin{align*} 
		&\left\lvert R_1\right\rvert \leq\left\lVert\varphi\right\rVert_{\mathscr{L}^p_\mathbb{R}}\left(\mathbb{E}\int_0^T\left\lvert\int_{t}^{\kappa_2(t)} \left\langle A_s\left(u^{n,m,l}(s)\right),z\right\rangle \mathrm{d} s\right\rvert^q\mathrm{d} t\right)^{1/q}\\&\leq \left\lVert\varphi\right\rVert_{\mathscr{L}^p_\mathbb{R}}\left\lVert z\right\rVert_V\\&\times\left[\mathbb{E}\int_0^T\left\lvert\int_{0}^{T} \left\lVert A_s\left(u^{n,m,l}(s)\right)\right\rVert_{V^*}^q\lambda(s)^{1-q} \mathrm{d} s\right\rvert\cdot\left\lvert\int_{t}^{\kappa_2(t)}\lambda(s)\mathrm{d} s\right\rvert^{q/p}\mathrm{d} t\right]^{1/q}\\&\leq 
		\left\lVert\varphi\right\rVert_{\mathscr{L}^p_\mathbb{R}} \left\lVert z\right\rVert_VT^{1/q}\max_{1\leq i\leq m}\left(\int_{t_{i-1}}^{t_i}\lambda(s)\mathrm{d} s\right)^{1/p} \\&\quad  \times\left(\mathbb{E}\int_0^T \left\lVert A_s\left(u^{n,m,l}(s)\right)\right\rVert_{V^*}^q\lambda(s)^{1-q}\mathrm{d} s\right)^{1/q}\to 0 \, . 
	\end{align*} 
	For $R_2$, using Cauchy-Schwartz's inequality we have 
	\begin{align*} 
		R_2^2 &=\left\lvert \mathbb{E}\int_0^T\varphi(t)\int_{t}^{\kappa_2(t)} \left(\Pi_n\tilde{B}^m_s\left(u^{n,m,l}(\kappa_1(s))\right)\tilde{\Pi}_l\mathrm{d} W_s,z\right)_H\mathrm{d} t\right\rvert^2\\ 
		&\leq\left\lVert\varphi\right\rVert^2_{\mathscr{L}_\mathbb{R}^2}\mathbb{E}\int_0^T\left\lvert\int_{t}^{\kappa_2(t)}\left(\Pi_n\tilde{B}^m_s\left(u^{n,m,l}(\kappa_1(s))\right)\tilde{\Pi}_l\mathrm{d} W_s,z\right)_H\right\rvert^2\mathrm{d} t\\&\leq \left\lVert\varphi\right\rVert^2_{\mathscr{L}^2_\mathbb{R}}\left\lVert z\right\rVert^2_H\mathbb{E}\int_0^T\int_{t}^{\kappa_2(t)} \left\lVert \Pi_n\tilde{B}^m_s\left(u^{n,m,l}(\kappa_1(s))\right)\tilde{\Pi}_l\right\rVert_{2}^2 \mathrm{d} s\,\mathrm{d} t\\&\leq 
		\delta_m\left\lVert\varphi\right\rVert^2_{\mathscr{L}^2_\mathbb{R}}\left\lVert z\right\rVert^2_H \mathbb{E}\int_0^T \left\lVert \Pi_n\tilde{B}^m_t\left(u^{n,m,l}(\kappa_1(t))\right)\tilde{\Pi}_l\right\rVert_{2}^2\mathrm{d} t, 
	\end{align*} 
since $\left\lVert\Pi_n\tilde{B}^m_s\left(u^{n,m,l}(\kappa_1(s))\right)\tilde{\Pi}_l\right\rVert_{2}^2$ is constant for $s\in (t,\kappa_2(t))$. Thus by Step 1, $R_2\to 0$ when $m\to \infty$. Finally the computation for $R_3$ is as follows: 
		\begin{align*} 
			&R_3^2 = \left\lvert\mathbb{E}\int_0^T\varphi(t)\int_{t}^{\kappa_2(t)}\int_E \left(\Pi_n\tilde{F}^{m,l}_s\left(u^{n,m,l}(\kappa_1(s)),\xi\right),z\right)_H\tilde{N}(\mathrm{d} s,\mathrm{d}\xi)\mathrm{d} t\right\rvert^2\\ 
			&\leq \left\lVert\varphi\right\rVert_{\mathscr{L}_\mathbb{R}^2}^2\mathbb{E}\int_0^T\left\lvert \int_{t}^{\kappa_2(t)}\int_E \left(\Pi_n\tilde{F}^{m,l}_s\left(u^{n,m,l}(\kappa_1(s)),\xi\right),z\right)_H\tilde{N}(\mathrm{d} s,\mathrm{d}\xi)\right\rvert^2\mathrm{d} t\\&\leq	\left\lVert\varphi\right\rVert_{\mathscr{L}_\mathbb{R}^2}^2\left\lVert z\right\rVert^2_H\mathbb{E}\int_0^T \int_{t}^{\kappa_2(t)}\int_E \left\lVert\Pi_n\tilde{F}^{m,l}_s\left(u^{n,m,l}(\kappa_1(s)),\xi\right)\right\rVert^2_H\nu(\mathrm{d}\xi)\mathrm{d} s\,\mathrm{d} t\\&\leq 	\delta_m\left\lVert\varphi\right\rVert_{\mathscr{L}_\mathbb{R}^2}^2\left\lVert z\right\rVert^2_H\mathbb{E}\int_0^T\int_E \left\lVert \Pi_n\tilde{F}^{m,l}_t\left(u^{n,m,l}(\kappa_1(t)),\xi\right)\right\rVert_{H}^2\nu(\mathrm{d}\xi)\mathrm{d} t
		\end{align*} 
In the last inequality, we used the fact that $\left\lVert\Pi_n\tilde{F}^{m,l}_s\left(u^{n,m,l}(\kappa_1(s)),\xi\right)\right\rVert^2_H$ is constant for $s\in (t,\kappa_2(t))$. So by Step 1,  $R_3\to 0$ when $m\to \infty$. 
	Now we have proven that the limit of right hand side of equation \eqref{eq:Ji Ri's imp} is the right hand side of equation \eqref{effect of v inftyimp}. By the fact that $u^{n,m,l}\rightharpoonup \bar{u}_\infty$, we deduce the similar result for the left hand side, 
	so equation \eqref{effect of v inftyimp} is proved. 
	It remains to prove \eqref{u inftyimp}. Both sides of this equation belong to the space $L^1(\Omega;\mathbb{R})$, 
	 So it is sufficient to prove that for every $\psi \in L^\infty(\Omega;\mathbb{R})$ and for $z\in V_N$ 
	\begin{align*} 
		\mathbb{E}\left[\psi\left(u_{T\infty},z\right)_H\right]&=\mathbb{E}\left[\psi(\zeta ,z)_H\right]+\mathbb{E}\left[\psi\int_0^T\left\langle 
		a_\infty(s),z\right\rangle \mathrm{d} s\right]\\&\quad+\mathbb{E}\left[\psi\int_0^T\left(b_\infty(s)\mathrm{d} W_s,z\right)_H\right]\\&\quad +\mathbb{E}\left[\psi\int_0^T\int_E\left(f_\infty(s,\xi),z\right)_H\tilde{N}(\mathrm{d} s,\mathrm{d}\xi)\right]. 
	\end{align*} 
	Fix $N\in\mathbb{N}$ and $z\in V_N$. Using equation \eqref{int form of imp}, we get for $n\geq N$ that 
	\begin{align*} 
		(u^{n,m,l}(T),z)_H&=(\zeta ,z)_H\mathbf{1}_{\left\lbrace t>t_0\right\rbrace}+\int_0^{T}\left\langle A_s\left(u^{n,m,l}(s)\right),z\right\rangle_H\mathrm{d} s\\&\quad +\int_0^{T}\left(z,\Pi_n\tilde{B}^m_s\left(u^{n,m,l}\left(\kappa_1(s)\right)\right)\tilde{\Pi}_l\mathrm{d} W_s\right)_H\\ 
		&\quad+\int_0^{T}\int_E\left(\Pi_n\tilde{F}^{m,l}_s\left(u^{n,m,l}\left(\kappa_1(s)\right)\right),z\right)_H\tilde{N}(\mathrm{d} s,\mathrm{d}\xi) \, . 
	\end{align*} 
	Taking the inner products of both sides and $\psi$, we get for $n\geq N$ 
	\begin{equation*}
		\mathbb{E}\left[\psi(u^{n,m,l}(T),z)_H\right]=\mathbb{E}\left[\psi(\zeta ,z)_H\right]+\tilde{J}_1+\tilde{J}_2+\tilde{J}_3 \, , 
	\end{equation*} 
	where 
	\begin{align*} 
		\tilde{J}_1&=\mathbb{E}\left[\psi\int_0^T \left\langle A_s\left(u^{n,m,l}(s)\right),z\right\rangle \mathrm{d} s\right] \to \mathbb{E}\left[\psi\int_0^T \left\langle a_\infty(s),z\right\rangle \mathrm{d} s\right],\\ 
		\tilde{J}_2&=\mathbb{E}\left[\psi\int_0^T\left( \Pi_n\tilde{B}^m_s\left(u^{n,m,l}(\kappa_1(s))\right)\tilde{\Pi}_l\mathrm{d} W_s,z\right)_H\right]\\&\to \mathbb{E}\left[\psi\int_0^T\left( b_\infty(s)\mathrm{d} W_s,z\right)_H\right],\\ 
		\tilde{J}_3&=\mathbb{E}\left[\psi\int_0^T\int_E \left(\Pi_n\tilde{F}^{m,l}_s\left(u^{n,m,l}(\kappa_1(s)),\xi\right),z\right)_H\tilde{N}(\mathrm{d} s,\mathrm{d}\xi)\right]\\&\to\mathbb{E}\left[\psi\int_0^T\int_E \left(f_\infty(s,\xi),z\right)_H\tilde{N}(\mathrm{d} s,\mathrm{d}\xi)\right]. 
	\end{align*} 
	The limits are obtained similarly as the limits of the corresponding terms in \eqref{eq:Ji Ri's tilde}. So \eqref{u inftyimp} is proved. 
	
	\paragraph*{\textbf{Step 3. }}
	\textit{ Let $(n,m,l)$ be a subsequence which satisfies items \ref{converg of implicit:i}-\ref{converg of implicit:vi} of the previous step. Then for all $y\in \mathscr{G}_{\mathcal{BF}}$ 
		\[I_y(\zeta ,a_\infty,b_\infty,f_\infty)+\liminf \mathbb{E}\left\lVert u^{n,m,l}(T)\right\rVert_H^2-\mathbb{E}\left\lVert u_{T\infty}\right\rVert_H^2\leq 0.\]}
	\paragraph*{\text{Proof of Step 3. }} 
Let $\mathcal{T}_m$ and $\mathcal{S}_l$ be the operators defined in \eqref{Tm-define} and \eqref{Sl-define} respectively. Then
		\begin{equation*}
			\begin{split}
				\mathcal{T}_m(\Pi_nB_\cdot\left(u_{n,m,l}(\cdot)\right)\tilde{\Pi}_l)&=\Pi_n\tilde{B}^m_\cdot\left(u_{n,m,l}(\kappa_1(\cdot))\right)\tilde{\Pi}_l,
				\\\mathcal{S}_l\circ \mathcal{T}_m(\Pi_nF_\cdot\left(u_{n,m,l}(\cdot),\star\right))&=\Pi_n\tilde{F}^{m,l}_\cdot\left(u_{n,m,l}(\kappa_1(\cdot)), \star\right).
			\end{split}
		\end{equation*}
	 Define 
	\begin{align*} 
		&I^{n,m,l}(y):=\mathbb{E}\int_0^T \Bigg[\left\langle A_s\left(u^{n,m,l}(s)\right)-A_s(y_s), u^{n,m,l}(s)-y_s\right\rangle\\&\quad+\left\lVert \Pi_n\tilde{B}^m_s\left(u^{n,m,l}(\kappa_1(s))\right)\tilde{\Pi}_l-\mathcal{T}_m\left(\Pi_nB_\cdot(y_\cdot)\right)(s)\right\rVert_2^2\\& \quad+\int_E\left\lVert \Pi_n\tilde{F}^{m,l}_s\left(u^{n,m,l}(\kappa_1(s)),\xi\right)-\mathcal{S}_l\circ \mathcal{T}_m\left(\Pi_nF_\cdot(y_\cdot,\star)\right)(s,\xi)\right\rVert_H^2\nu(\mathrm{d} \xi)\Bigg]\mathrm{d} s. 
	\end{align*} 
Since the operator norms of $\mathcal{T}_m$ and $\mathcal{S}_l$ equal one, we get
		\begin{align*} 
			I^{n,m,l}(y)&\leq\mathbb{E}\int_0^T \Bigg[\left\langle A_s\left(u^{n,m,l}(s)\right)-A_s(y_s), u^{n,m,l}(s)-y_s\right\rangle\\&\quad+\left\lVert \Pi_n B_s\left(u^{n,m,l}(s)\right)\tilde{\Pi}_l-\Pi_nB_s(y_s)\right\rVert_2^2\\& \quad+\int_E\left\lVert \Pi_nF_s\left(u^{n,m,l}(s),\xi\right)-\Pi_nF_s(y_s,\xi)\right\rVert_H^2\nu(\mathrm{d} \xi)\Bigg]\mathrm{d} s \, . 
	\end{align*} 
	The monotonicity condition \ref{C1} implies that $I^{n,m,l}(y)\leq 0$. Summing up inequality \eqref{ghablan(4)imp} with respect to $i$, $1\leq i\leq m$, we get 
	\begin{align*} 
		&\mathbb{E}\left\lVert u^{n,m,l}(T)\right\rVert^2_H \leq \mathbb{E}\left\lVert \zeta \right\rVert_H^2 +\mathbb{E}\int_0^T\Bigg[ 2\left\langle u^{n,m,l}(s),A_s\left(u^{n,m,l}(s)\right)\right\rangle\\& +\left\lVert \Pi_n \tilde{B}^m_s\left(u^{n,m,l}(\kappa_1(s))\right)\tilde{\Pi}_l\right\rVert_2^2+\int_E\left\lVert \Pi_n \tilde{F}^{m,l}_s\left(u^{n,m,l}(\kappa_1(s)),\xi\right)\right\rVert_H^2\nu(\mathrm{d} \xi)\Bigg]\mathrm{d} s. 
	\end{align*} 
	Thus 
	\begin{align}\label{I^n,m,l} 
		0\geq I^{n,m,l}(y)&\geq \mathbb{E}\left\lVert u^{n,m,l}(T)\right\rVert_H^2-\mathbb{E}\left\lVert \zeta \right\rVert_H^2+\mathbb{E}\int_0^T 2\left\langle A_s(y_s),y_s\right\rangle\mathrm{d} s\nonumber\\&\quad+J_1+J_2-2(J_3+J_4+J_5+J_6) \, , 
	\end{align} 
	where 
	\begin{align*} 
		J_1&:=\mathbb{E}\int_0^T\left\lVert \mathcal{T}_m\left(\Pi_n B_\cdot(y_\cdot)\right)(s)\right\rVert_2^2\mathrm{d} s\to \mathbb{E}\int_0^T\left\lVert B_s(y_s)\right\rVert_2^2\mathrm{d} s, 
		\\J_2&:=\mathbb{E}\int_0^T\int_E\left\lVert \mathcal{S}_l\circ \mathcal{T}_m\left(\Pi_n F_\cdot(y_\cdot,\star)\right)(s,\xi)\right\rVert_H^2\nu(\mathrm{d} \xi)\mathrm{d} s\\&\to \mathbb{E}\int_0^T\int_E\left\lVert F_s(y_s,\xi)\right\rVert_H^2\nu(\mathrm{d} \xi)\mathrm{d} s,\\ 
		J_3&:=\mathbb{E}\int_0^T\left\langle A_s\left(u^{n,m,l}(s)\right),y_s\right\rangle\mathrm{d} s\to \mathbb{E}\int_0^T \left\langle a_\infty(s),y_s\right\rangle\mathrm{d} s,\\ 
		J_4&:=\mathbb{E}\int_0^T\left\langle A_s(y_s),u^{n,m,l}(s)\right\rangle\mathrm{d} s\to \mathbb{E}\int_0^T\left\langle A_s(y_s),\bar{u}_\infty(s)\right\rangle\mathrm{d} s,\\ 
		J_5&:=\mathbb{E}\int_0^T\left\langle \mathcal{T}_m\left(\Pi_n B_\cdot(y_\cdot)\right)(s),\Pi_n\tilde{B}^m_s\left(u^{n,m,l}(\kappa_1(s))\right)\tilde{\Pi}_l\right\rangle_2\mathrm{d} s\\&\to \mathbb{E}\int_0^T\int_E\left\langle B_s(y_s),b_\infty(s)\right\rangle_2\mathrm{d} s ,\\ 
		J_6&:=\mathbb{E}\int_0^T\int_E\Big( \mathcal{S}_l\circ \mathcal{T}_m\left(\Pi_n F_\cdot(y_\cdot,\star)\right)(s,\xi),\\&\qquad \qquad\qquad\Pi_n \tilde{F}^{m,l}_s\left(u^{n,m,l}(\kappa_1(s)),\xi\right)\Big)_H\nu(\mathrm{d} \xi)\mathrm{d} s\\&\to \mathbb{E}\int_0^T\int_E\left( F_s(y_s,\xi),f_\infty\left(s,\xi\right)\right)_H\nu(\mathrm{d} \xi)\mathrm{d} s. 
	\end{align*} 
	Since $u^{n,m,l}(T)\rightharpoonup u_{T\infty}$, we can write 
	\begin{equation}\label{d-implicit}
	d:=\liminf \mathbb{E}\left\lVert u^{n,m,l}(T)\right\rVert_H^2-\mathbb{E}\left\lVert u_{T\infty}\right\rVert_H^2\geq 0.
	\end{equation}
	By Step 2 and Theorem \ref{theorem1of1}, $\bar{u}_\infty$ has an adapted c\`adl\`ag $H$-valued modification, denoted with $u_\infty$, such that almost surely for all $t\in[0,T]$ 
	\[u_\infty(t)=\zeta +\int_0^t a_\infty(s)\mathrm{d} s+\int_0^t b_\infty (s)\mathrm{d} W_s+\int_0^t\int_E f_\infty(s,\xi)\tilde{N}(\mathrm{d} s, \mathrm{d} \xi).\] 
	This equation for $t=T$ together with the equation \eqref{u inftyimp} implies $u_\infty(T)=u_{T\infty}$ a.s.. 
	By It\^o's formula (see Theorem \ref{theorem1of1}), we have 
	\begin{equation} \label{uTinfty-implicit}
	\begin{split}
		&\mathbb{E}\left\lVert u_{T\infty}\right\rVert_H^2 =\mathbb{E}\left\lVert u_\infty(T)\right\rVert_H^2\\& =\mathbb{E}\left\lVert \zeta \right\rVert_H^2+\mathbb{E}\int_0^T\left[ 2\left\langle a_\infty(s),\bar{u}_\infty(s)\right\rangle+\left\lVert b_\infty(s)\right\rVert_2^2+\int_E\left\lVert f_\infty(s,\xi)\right\rVert_H^2\nu(\mathrm{d} \xi)\right]\mathrm{d} s. 
		\end{split}
	\end{equation} 
	Therefore, by taking limit inferior of inequality \eqref{I^n,m,l} and using \eqref{d-implicit} and \eqref{uTinfty-implicit}, we conclude that
		\begin{align*} 
			0&\geq d +\mathbb{E}\int_0^T\left[ 2\left\langle a_\infty(s),\bar{u}_\infty(s)\right\rangle+\left\lVert b_\infty(s)\right\rVert_2^2+\int_E\left\lVert f_\infty(s,\xi)\right\rVert_H^2\nu(\mathrm{d} \xi)\right]\mathrm{d} s\\&+\mathbb{E}\int_0^T \left[2\left\langle A_s(y_s),y_s\right\rangle+\left\lVert B_s(y_s)\right\rVert_2^2+\int_E\left\lVert F_s(y_s,\xi)\right\rVert_H^2\nu(\mathrm{d} \xi)\right]\mathrm{d} s\\& -2\, \mathbb{E}\int_0^T \left[\left\langle a_\infty(s),y_s\right\rangle+\left\langle A_s(y_s),\bar{u}_\infty(s)\right\rangle\right]\mathrm{d} s \\&-2\, \mathbb{E}\int_0^T \left[\left\langle B_s(y_s),b_\infty(s)\right\rangle_2+\int_E\left( F_s(y_s,\xi),f_\infty\left(s,\xi\right)\right)_H\nu(\mathrm{d} \xi)\right]\mathrm{d} s \, , \\&=
			d+I_y(\zeta ,a_\infty,b_\infty,f_\infty). 
	\end{align*} 
	Hence by Theorem \ref{characterize}, $u_\infty$ must be a solution of equation \eqref{equation1} 
	which is unique. By setting $y=u_\infty$ one obtains that $d=0$. The weak convergence of 
	$u^{n,m,l}(T)$, combined with the fact that $d=0$, implies that a subsequence of $u^{n,m,l}(T)$ 
	converges strongly in $L^2(\Omega,\mathbb{P};H)$ to $u_\infty(T)$. Since there exist unique 
	limits for any convergent subsequences of $u^{n,m,l}$ and $u^{n,m,l}(T)$, and on the other hand 
	each subsequence of these sequences has a convergent subsequence, the whole sequences $u^{n,m,l}$ and $u^{n,m,l}(T)$ converge. Thus the proof of theorem is complete. 
\end{proof} 
\section{Example: Stochastic parabolic equations}\label{example}
In this section, we illustrate an example for the provided set-up and we conduct a simulation using both the explicit and implicit schemes. We consider the equation
\begin{equation} \label{eq:example}
	\begin{split}\mathrm{d} u_t&=\left(\frac{1}{2}\partial_{xx} u_t-u_t\right)\mathrm{d} t+\theta\partial_x u_t\mathrm{d} W_t\\&\quad +\int_{[1,\infty)^2} \frac{u_t }{\xi_1\xi_2}\tilde{N}(\mathrm{d} t, \mathrm{d} \xi_1, \mathrm{d} \xi_2)\end{split}
\end{equation}
with given some initial condition $u_0(x)$ for $x\in[0,1]$.

To establish the provided framework, first we introduce the Gelfand triple $V=H^{1,2}_0(0,1)\subset H := L^2(0,1)\subset V^\star=H^{-1,2}(0,1)$ (see \cite{pardoux2007stochastic}). Here, $W$ is standard one-dimensional Brownian motion, $\theta\in(-1,1)$, and $\tilde{N}$ is compensated Poisson measure on $[0,T]\times[0,1]\times [1,\infty)^2$ with Lebesgue intensity measure. Now we define the corresponding operators
\[A:H^{1,2}_0(0,1)\to H^{-1,2}(0,1), \quad Au:= \frac{1}{2}u''-u,\]
\[B: H^{1,2}_0(0,1)\to L^2(0,1), \quad Bu:=\theta u',\]
\[F:H^{1,2}_0(0,1) \times [1,\infty)^2\to L^2(0,1), \quad F(u,\xi_1, \xi_2):=\frac{u}{\xi_1\xi_2}. \]

Let $V_n$ be the subspace of $V$ consisting of continuous piecewise linear functions like $v$ which are affine linear on $[(k-1)/n, k/n]$ for all $k\in \{1,2,3, \ldots, n\}$ and $v(0)=v(1)=0$. 
Consider the basis $\{e_1,\ldots, e_{n-1}\} $ for $V_n$, where
\[ e_k(x):=\mathbf{1}_{\left[\frac{k-1}{n},\frac{k+1}{n}\right]}(x) \left[1-|nx-k| \right].\]
So for $v\in V_n$, $[v]:=[v(1/n), \ldots, v((n-1)/n)]$ is the representation of $v$ with respect to the basis $\{e_1,\ldots, e_{n-1}\} $.

We have
\[(e_i,e_j)_H=\begin{dcases}
	\frac{2}{3\times n}& i=j\\
	\frac{1}{6\times n}& |i-j|=1\\
	0 & |i-j|>1
\end{dcases}
\qquad(e^\prime_i,e^\prime_j)_H=\begin{dcases}
	2n& i=j\\
	-n& |i-j|=1\\
	0 & |i-j|>1
\end{dcases}
\]
To represent this inner product in matrix form, we define a matrix $M=T\left(\frac{2}{3}, \frac{1}{6}, \frac{1}{6}\right)$ where
\[T(a,b,c) := \begin{pmatrix}
	a\,\, & b \\
	c\,\, & a & b \\
	& c & \ddots & \ddots \\
	& & \ddots & \ddots & b \\
	& & & c & a
\end{pmatrix}\]
is tridiagonal ${(n-1)\times (n-1)}$ matrix.
Thus, for $u,v\in V_n$
\[(u,v)_H=\frac{1}{n}[u]M[v]^T\]

The columns of $\sqrt{nM^{-1}}[e_1,\ldots , e_{n-1}]^T$ consist an orthonormal basis for $V_n$ with respect to inner product in $H$. So
\[\sum_{k=1}^{n-1}\|e_k\|_V^2=n^2\mathbf{tr}\left(M^{-1}T(2,-1,-1)\right)\simeq Cn^3\]
and the condition for convergence in \cite[Theorem 2.8]{gyongy2005discretization} reads as $n^3/m\to 0$. 

We need to obtain $\Pi_nv$, $\Pi_n v''$, $\Pi_nv'$.
For $v\in V^\star$ and $k\in \{1,2, 3, \ldots, n-1\}$, we have
\[\begin{split}\frac{1}{n}[\Pi_nv]M[e_k]^T&=\left( \Pi_n v, e_k\right)_H=\left\langle v, e_k\right\rangle\\&=\int_{(k-1)/n}^{(k+1)/n}v(r)e_k(r)\mathrm{d} r\end{split}\]
Since $[e_k]_i=\mathbf{1}_{k=i}$, we get
\[[\Pi_nv]=n[\left\langle v, e_1\right\rangle, \left\langle v, e_2\right\rangle, \ldots, \left\langle v, e_{n-1}\right\rangle]M^{-1}\]
For $v\in V_n$, we have
\[\begin{split}\left\langle v', e_k\right\rangle&=\int_0^1 v'(r)e_k(r)\mathrm{d} r=\frac{v((k+1)/n)-v((k-1)/n)}{2}.\end{split}\]
and
\[\begin{split}\left\langle v'', e_k\right\rangle&=\int_0^1 v''(r)e_k(r)\mathrm{d} r=-\int_0^1v'(r)e^\prime_k(r)\mathrm{d} r\\&-n\int_{(k-1)/n}^{k/n}v'(r)\mathrm{d} r+n\int_{k/n}^{(k+1)/n}v'(r) \mathrm{d} r\\&=n\left[v((k-1)/n)-2v(k/n)+v((k+1)/n)\right]
\end{split}\]
So for $v\in V_n$, we obtain the matrix representations of $\Pi_nv'$ and $\Pi_nv''$ as follows,
\[[\Pi_nv']=\frac{n^2}{2} [v]T\left(0,-1, 1\right) M^{-1}=:[v]M'.\]
and
\[[\Pi_nv''] = n^2[v]T(-2,1,1 )M^{-1} =: [v]M''\]
Let us estimate $\mathcal{C}(n)$. We have for $v\in V_n$,
\[\begin{split}\|\Pi_n A(v) \|_H^2&=\frac{1}{n}\left(\frac{1}{2}[v]M''-[v]\right)M\left(\frac{1}{2}M''[v]^T-[v]^T\right)\\&\leq C n^3[v]T(-2,1,1 )M^{-1}T(-2,1,1 )[v]^T\end{split}\]
and for $u\in V_n$
\[\left\langle v'',u\right\rangle=\sum_{k=1}^{n-1}n\left[v((k-1)/n)-2v(k/n)+v((k+1)/n)\right]u(k/n).\]
For $u\in V_n$ with $[u]=[v]T(-2,1,1)$, we have
\[\left\langle A(v),u\right\rangle= \left\langle \frac{1}{2}v''-v,u\right\rangle=\frac{n}{2} [u][u]^T+[v]T(2,-1,-1)[v]^T\geq \frac{n}{2} [u][u]^T, \]
\[ \|u\|_V^2=n[u]T(2, -1,-1)[u]^T\leq 4 n[u][u]^T. \]
So \[\begin{split}\|A(v)\|^2_{V^\star}&\geq \frac{|\left\langle A(v),u\right\rangle|^2}{\|u\|_V^2}\geq \frac{n}{4}[u][u]^T\\&\geq \frac{\|M\|n^3[u]M^{-1}[u]^T}{4n^2} \\&= \frac{\|M\|}{4Cn^2}\|\Pi_n A(v) \|_H^2\end{split}\]
Therefore, $\mathcal{C}(n)\leq 4Cn^2/\|M\|$ and the condition in Theorem \ref{explicitconvergence} reads as $n^2/m\to 0$ which is significantly weaker than the condition in \cite[Theorem 2.8]{gyongy2005discretization}.
\begin{figure}[h]
	\begin{center}
		\includegraphics[scale=.44]{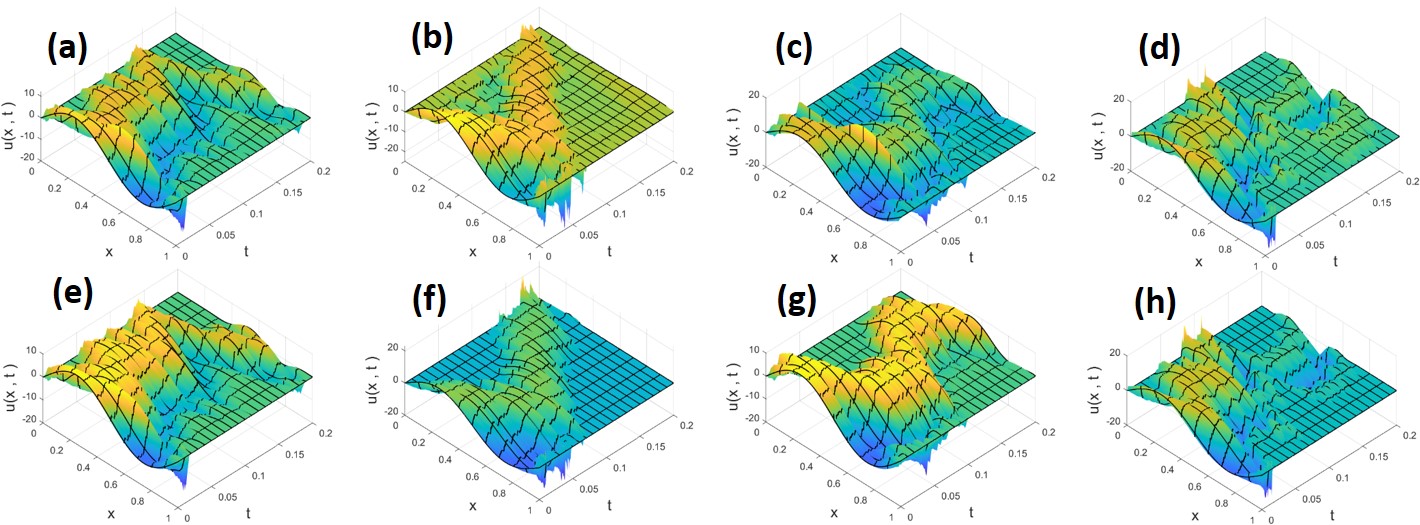}
	\end{center}
	\caption{Simulated solution of equation~\eqref{eq:example} for (a) explicit scheme with $n=50$ and $m=2500$ (b) explicit scheme with $n=200$ and $m=40000$ (c) implicit scheme with $n=50$ and $m=2500$ (d) implicit scheme with $n=200$ and $m=40000$. Plots (e) to (h) show the simulated solution without the Poisson noise (only Brownian noise) with the same parameters. We have set the initial condition $u_0(x) = 10\sin(2\pi x)$ and $\theta=.5$.In all cases we set $l=10$.}
	\label{fig:simulation}
\end{figure}

We also set $E^l := [1,l]\times[1,l]$ and partition $E^l$ by dividing $[1,l]$ into equal subintervals of length $\frac{1}{l}$. 

So explicit scheme becomes
\[[u^n_{m,l}(t_0)]=[\Pi_n u_0(x)],\]
\[[u^n_{m,l}(t_{k+1})]= [u^n_{m,l}(t_{k})] \left(I + \delta_m\left(\frac{1}{2}M''-I\right) + \theta \sqrt{\delta_m} \mathcal{N}M' + \tilde{Z}^l I \right) \]
where $\mathcal{N}$ is standard normal random variable and 
\[ \tilde{Z}^l = \sum_{j=1}^{r_l} \tilde{F}^l_j . (Z^l_j-\delta_m |E^l_j|)\]
and $Z^l_j$ is a Poisson random variable with parameter $\delta_m |E^l_j|$ and $\tilde{F}^l_j = \int_{E^l_j} \frac{d\xi_1d\xi_2}{\xi_1\xi_2}$.

For the implicit scheme, in each step we have to solve the following equation
\[[u^n_{m,l}(t_{k+1})]= [u^n_{m,l}(t_{k+1})] \delta_m\left(\frac{1}{2}M''-I\right) + [u^n_{m,l}(t_{k})] \left(I + \theta \sqrt{\delta_m} \mathcal{N}M' + \tilde{Z}^l I \right) \]
which leads to
\[[u^n_{m,l}(t_{k+1})]= [u^n_{m,l}(t_{k})] \left(I + \theta \sqrt{\delta_m} \mathcal{N}M' + \tilde{Z}^l I \right) \left( I - \delta_m\left(\frac{1}{2}M''-I\right)\right)^{-1} \]

Figure~\ref{fig:simulation} shows the graph of the solutions resulting from the explicit and implicit schemes for various values of parameters for both with and without the Poissonian noise.

\subsection*{Acknowledgements} 
We wish to thank Prof. Dr. Wilhelm Stannat, Professor of Mathematics at Technical University of Berlin, for assistance with editing, and for his useful comments that greatly improved the manuscript.

\bibliographystyle{plain} 
\bibliography{biblio} 
\end{document}